\documentclass{amsart}

\usepackage{a4wide,amssymb,amsrefs,lscape,mathrsfs,times}

\usepackage[arrow,matrix,curve]{xy}\SilentMatrices

\def\brk#1{\left\langle#1\right\rangle}

\def\1{^{-1}}
\def\cref#1#2#3{\left(#2\right.\left|\ #3\right)_{#1}}

\def\fp{f_+}
\def\fm{f_\cdot}
\def\php{\ph_+}
\def\phm{\ph_\cdot}
\def\gammap{\gamma_+}
\def\gammam{\gamma_\cdot}
\def\sigmap{\sigma_+}
\def\sigmam{\sigma_\cdot}

\def\sbar#1{\bar{#1\mathstrut}}
\def\ol#1{\overline{#1\mathstrut}}

\def\shat#1{\hat{#1\mathstrut}}
\def\wh#1{\widehat{#1\mathstrut}}

\def\stil#1{\tilde{#1\mathstrut}}
\def\wt#1{\widetilde{#1\mathstrut}}

\def\rarr#1{{\begin{array}{r}#1\end{array}}}

\def\smat#1{\begin{smallmatrix}#1\end{smallmatrix}}

\def\aass#1#2#3{\left\langle{#1},{#2},{#3}\right\rangle}
\def\acom#1#2{\left\{{#1},{#2}\right\}}
\def\coma#1#2#3#4{\left\langle\smat{#1&#2\\#3&#4}\right\rangle}
\def\mass#1#2#3{\left[{#1},{#2},{#3}\right]}
\def\ldis#1#2#3{\left[{#1}\smat{{#2}\\{#3}}\right\rangle}
\def\rdis#1#2#3{\left\langle\smat{{#1}\\{#2}}{#3}\right]}

\let\x\times

\let\ph\varphi
\let\le\leqslant

\let\lim\varprojlim

\let\then\Rightarrow

\let\comm\circlearrowleft

\makeatletter

\let\c@equation\c@subsection

\def\thmhead#1#2#3{%
 {\the\thm@notefont(\thmnumber{
    \@upn{#2}})}
   \thmname{#1}%
   \thmnote{ {\the\thm@notefont(#3)}}}

\makeatother

\def\theodef#1{\newtheorem{#1}[subsection]{#1}}

\theodef{Lemma}
\theodef{Corollary}
\theodef{Proposition}
\theodef{Theorem}

\theoremstyle{definition}

\theodef{Definition}

\theodef{Remark}

\theodef{Example}
\theodef{Examples}

\def\xto#1{\xrightarrow[]{#1}}

\def\fb{{\boldsymbol f}}

\def\sA{{\mathscr A}}
\def\sR{{\mathscr R}}

\begin{document}

\title{Third Mac Lane cohomology via categorical rings}

\author{M. Jibladze}

\author{T. Pirashvili}

\maketitle

\section*{Introduction}

In the fifties, Saunders Mac~Lane invented a cohomology theory of
rings using the cubical construction introduced earlier by Eilenberg
and himself to calculate stable homology of Eilenberg-Mac~Lane
spaces. As shown in \cite{PW}, this theory coincides with the
topological Hochschild cohomology for Eilenberg-Mac~Lane ring
spectra. In particular, the third dimensional cohomology group is
expected to provide classification of 2-types of ring spectra. Some
algebraic models for such 2-types have been constructed in
\cite{BJP}. In this paper we consider one such algebraic model of
different kind which in our opinion is especially straightforwardly
related to 3-cocycles in Mac~Lane cohomology.

This is the notion of categorical ring---a category carrying the
structure of a ring up to some natural isomorphisms satisfying
certain coherence conditions. Our axioms for the categorical ring
present a slightly modified version of the notion of Ann-category
due to Quang \cite{queng}. Axioms we use reflect defining relations
of Mac~Lane 3-cocycles. Our main result is Theorem \ref{main} which
asserts that for any ring $R$ and any $R$-bimodule $B$ there is a
bijection
\[
H^3(R;B)\approx\mathrm{Crext}(R;B)
\]
between the third Mac~Lane cohomology group of $R$ with coefficients
in $B$ and equivalence classes of categorical rings $\sR$ with
$\pi_0(\sR)=R$, $\pi_1(\sR)=B$, and the induced bimodule structure
coinciding with the original one.

In \cite{queng}, Ann-categories of a particular kind---the so called
regular ones---are considered. These correspond to the ring spectra
whose underlying spectrum splits into a product of
Eilenberg-Mac~Lane spectra. It is shown in \cite{queng} that regular
Ann-categories are classified by the third Shukla cohomology group
\cite{shukla}. The latter is the Barr-Beck-Quillen cohomology group
for the category of associative rings \cite{BPA}.

Difference between the above two cases is quite subtle. In fact,
Shukla and Mac~Lane cohomologies are isomorphic up to dimension 2
(in dimensions 0 and 1 both also coincide with the Hochschild
cohomology; this coincidence extends to dimension 2 if the
underlying abelian group of $R$ is free). The third Shukla
cohomology group embeds into the third Mac~Lane cohomology group,
and failure of isomorphism is measured by an explicit obstruction
furnished by the following exact sequence (see \cites{ill,louvain},
\cite{JP}, \cite{BPA}):
\[
0\to H^3_{\mathrm{Shukla}}(R;B)\to H^3(R;B)\to H^0(R;{}_2B).
\]
It can be shown that in terms of categorical rings, the above map
$H^3(R;B)\to H^0(R;{}_2B)$ sends an element of $H^3(R;B)$
represented by a categorical ring $\sR$ to the element
\[
0\to0+0\xto{\acom00}0+0\to0
\]
of Aut$_\sR(0)=B$, where $0$ is the neutral object with respect to
the additive structure, $\acom xy:x+y\to y+x$ is the commutativity
constraint, and the isomorphisms between $0$ and $0+0$ are the
canonical ones.

\section{Recollections on symmetric categorical groups}\label{catg}

Let us begin by recalling

\begin{Definition}
A \emph{categorical group} $\sA$ is a groupoid equipped with a
monoidal structure---i.~e. a bifunctor $+:\sA\x\sA\to\sA$, an object
$0\in\sA$ and natural isomorphisms $\aass abc:(a+b)+c\to a+(b+c)$,
$\lambda(a):0+a\to a$ and $\rho(a):a+0\to a$ satisfying the Mac~Lane
coherence conditions---together with a choice, for each object
$a\in\sA$, of another object $-a\in\sA$ and of an isomorphism
$\iota(a):-a+a\to0$.
\end{Definition}

Categorical groups are also known in the literature under the name
of \emph{Picard categories}.

We will need some specific facts and auxiliary notation concerning
symmetric categorical groups.

It is easy to see that for any monoidal functor
$\fb=(f,\fp,f_0):\sA\to\sA'$ to a categorical group, the canonical
isomorphism $f_0:f(0)\to0$ is determined by the rest of the
structure. Namely, $f_0$ is equal to the composite
\begin{equation}\label{zero}
\begin{aligned}
\xymatrix@!C=6ex{%
f(0)%
\ar[d]_-{\lambda(f(0))}%
&&0%
\\
0+f(0)%
\ar[d]_-{\iota(f(0))\1+f(0)}%
&&-f(0)+f(0)%
\ar[u]_-{\iota(f(0))}%
\\
(-f(0)+f(0))+f(0)%
\ar[dr]_-{\aass{-f(0)}{f(0)}{f(0)}\qquad}%
&&-f(0)+f(0+0)%
\ar[u]_-{-f(0)+f(\lambda(0))}%
\\
&-f(0)+(f(0)+f(0))%
\ar[ur]_-{\qquad-f(0)+\fp(0,0)}%
}%
\end{aligned}
\end{equation}

We will use the well-known fact that a monoidal functor
$\fb:\sA\to\sA'$ between categorical groups is an equivalence if and
only if it induces an isomorphism on $\pi_0$ and $\pi_1$. Here,
$\pi_1$ of a categorical group is defined to be the automorphism
group of its neutral object, and the homomorphism
$f_\#:\pi_1(\sA)\to\pi_1(\sA')$ induced by a monoidal functor $\fb$
assigns to $\beta:0_\sA\to0_\sA$ the composite
\[
\xymatrix@1{%
0_{\sA'}\ar[r]^-{f_0\1}&f(0_\sA)\ar[r]^{f(\beta)}&f(0_\sA)\ar[r]^-{f_0}&0_{\sA'}.%
}%
\]

It is equally well known that in a categorical group $\sA$,
hom$(x,y)$ has a structure of a bitorsor under $\pi_1(\sA)$ for any
isomorphic objects $x$, $y$ of $\sA$. In particular, for any two
parallel morphisms $\alpha,\alpha':x\to y$ there exists a
\emph{unique} $\beta\in\pi_1(\sA)$ making the diagram
\[
\xymatrix{%
x+0%
\ar[r]^{\alpha+\beta}%
\ar[d]_{\rho(x)}%
&y+0%
\ar[d]^{\rho(y)}%
\\
x%
\ar[r]^{\alpha'}
&y%
}%
\]
commute. It will be more convenient for us to depict such
circumstances by a diagram of the form
\[
\xymatrix@C=4em{%
x%
\ar@/^1pc/[r]^{\alpha}%
\ar@{}[r]|{\beta}%
\ar@/_1pc/[r]_{\alpha'}%
&y.%
}%
\]

Note that exchanging order of $\alpha$ and $\alpha'$ introduces a
sign, i.~e. one has
\[
\xymatrix@C=4em{%
x%
\ar@/^1pc/[r]^{\alpha}%
\ar@{}[r]|{\beta}%
\ar@/_1pc/[r]_{\alpha'}%
&y%
}%
\qquad\iff\qquad%
\xymatrix@C=4em{%
x%
\ar@/^1pc/[r]^{\alpha'}%
\ar@{}[r]|{-\beta}%
\ar@/_1pc/[r]_{\alpha}%
&y;%
}%
\]
in the diagrams that we will encounter, this order is not specified
as it can be unambiguously recovered from the context.

We will need the following simple fact concerning this formalism.

\begin{Proposition}\label{holef}
Let $\fb=(f,\fp):\sA\to\sA'$ be a monoidal functor between
categorical groups. Then for any parallel arrows
$\alpha,\alpha':x\to y$ in $\sA$ one has
\[
\xymatrix@C=4em{%
x%
\ar@/^1pc/[r]^{\alpha}%
\ar@{}[r]|{\beta}%
\ar@/_1pc/[r]_{\alpha'}%
&y%
}%
\qquad\then\qquad%
\xymatrix@C=4em{%
f(x)%
\ar@/^1pc/[r]^{f(\alpha)}%
\ar@{}[r]|{f_\#(\beta)}%
\ar@/_1pc/[r]_{f(\alpha')}%
&f(y).%
}%
\]
\end{Proposition}

\begin{proof}
There is a commutative diagram
\[
\xymatrix@C=4em{%
f(x)+0%
\ar[r]^{f(\alpha)+f_\#(\beta)}%
\ar[d]_{f(x)+f_0\1}%
&f(y)+0%
\ar[d]^{f(y)+f_0\1}%
\\
f(x)+f(0)%
\ar[d]_{\fp(x,0)}%
\ar[r]^{f(\alpha)+f(\beta)}%
&f(y)+f(0)%
\ar[d]^{\fp(y,0)}%
\\
f(x+0)%
\ar[r]^{f(\alpha+\beta)}%
\ar[d]_{f(\rho(x))}%
&f(y+0)%
\ar[d]^{f(\rho(y))}%
\\
f(x)%
\ar[r]^{f(\alpha')}
&f(y).%
}%
\]
By the aforementioned uniqueness, the proposition follows.
\end{proof}

\begin{Corollary}\label{holep}
For any parallel arrows $\alpha_i,\alpha'_i:x_i\to y_i$, $i=1,2$, in
a braided (in particular, symmetric) categorical group $\sA$ one has
\[
\xymatrix@C=4em{%
x_1%
\ar@/^1pc/[r]^{\alpha_1}%
\ar@{}[r]|{\beta_1}%
\ar@/_1pc/[r]_{\alpha'_1}%
&y_1%
},%
\quad%
\xymatrix@C=4em{%
x_2%
\ar@/^1pc/[r]^{\alpha_2}%
\ar@{}[r]|{\beta_2}%
\ar@/_1pc/[r]_{\alpha'_2}%
&y_2%
}%
\qquad\then\qquad
\xymatrix@C=4em{%
x_1+x_2%
\ar@/^1pc/[r]^{\alpha_1+\alpha_2}%
\ar@{}[r]|{\beta_1+\beta_2}%
\ar@/_1pc/[r]_{\alpha'_1+\alpha'_2}%
&y_1+y_2.%
}
\]
\end{Corollary}

\begin{proof}
It is well known that for a braided category $\sA$ the functor
$+:\sA\x\sA\to\sA$ acquires a monoidal structure (in fact it is
known that for any monoidal category there is a one-to-one
correspondence between braidings and monoidal functor structures on
$+$). Since obviously $+_\#(\beta_1,\beta_2)=\beta_1+\beta_2$ for
any $\beta_1,\beta_2\in\pi_1(\sA)$, the statement follows from
\eqref{holef}.
\end{proof}

For any four objects $a$, $b$, $c$, $d$ of a symmetric categorical
group $\sA$, by
\[
\coma abcd:(a+b)+(c+d)\to (a+c)+(b+d)
\]
will be denoted the composite canonical isomorphism in the
commutative diagram
\[
\xymatrix@!C=2em@!R=3em{
&&(a+b)+(c+d)\ar[drr]_{\aass ab{c+d}}\ar[dll]^{\aass{a+b}cd\1}\\
((a+b)+c)+d\ar[dr]^{\aass abc+d}&&&&a+(b+(c+d))\ar[dl]_{a+\aass bcd\1}\\
&(a+(b+c))+d\ar[rr]^{\aass a{b+c}d}\ar[d]_{(a+\acom
bc)+d}&&a+((b+c)+d)\ar[d]^{a+(\acom bc+d)}\\
&(a+(c+b))+d\ar[rr]_{\aass a{c+b}d}\ar[dl]^{\aass
acb\1+d}&&a+((c+b)+d)\ar[dr]_{a+\aass cbd}\\
((a+c)+b)+d\ar[rrd]^{\aass{a+c}bd}&&&&a+(c+(b+d))\ar[lld]_{\aass ac{b+d}\1}\\
&&(a+c)+(b+d) }
\]

\section{Categorical rings}

Our algebraic models for 2-types of ring spectra are certain
bimonoidal categories which we call categorical rings. They can be
called ``rings up to coherent isomorphisms'', in the sense that (a)
isomorphism classes of objects of a categorical ring form an
associative ring; and (b) the structure of a categorical ring on a
category is an equivalence invariant, i.~e. any equivalence between
a categorical ring and another category allows one to transfer the
categorical structure along it.

\begin{Definition}\label{2ring}
A \emph{categorical ring} is a symmetric categorical group $\sR$
together with a bifunctor $\sR\x\sR\to\sR$ (denoted by
juxtaposition), an object $1\in\sR$, and natural isomorphisms
\[
\mass rst:(rs)t\to r(st)
\]
(associativity),
\[
\lambda.(r):1r\to r,\ \ \rho.(r):r1\to r
\]
(left and right unitality),
\[
\ldis r{s_0}{s_1}:r(s_0+s_1)\to rs_0+rs_1
\]
(left distributivity),
\[
\rdis{r_0}{r_1}s:(r_0+r_1)s\to r_0s+r_1s
\]
(right distributivity).

It is required that the $\mass{}{}{}$ together with $\lambda.$ and
$\rho.$ constitute a monoidal structure (i.~e. the appropriate
pentagonal and triangular coherence diagrams commute for it) and
moreover the following diagrams commute for all possible objects of
$\sR$:
\[
\xymatrix@C=.5em@L=2ex
{
&r(s(t_0+t_1))\ar[rr]^{r\ldis s{t_0}{t_1}}&&r(st_0+st_1)\ar[dr]^{\ldis r{st_0}{st_1}}\\
(rs)(t_0+t_1)\ar[ur]^{\mass rs{t_0+t_1}}\ar[drr]_{\ldis{rs}{t_0}{t_1}}&&&&r(st_0)+r(st_1)\\
&&(rs)t_0+(rs)t_1\ar[urr]_{\mass rs{t_0}+\mass rs{t_1}}
}%
\]

\[
\xymatrix@C=3em@L=2ex
{
&(rs_0+rs_1)t\ar[dr]^{\rdis{rs_0}{rs_1}t}\\
(r(s_0+s_1))t\ar[d]_{\mass r{s_0+s_1}t}\ar[ur]^{\ldis
r{s_0}{s_1}t}&&(rs_0)t+(rs_1)t\ar[d]^{\mass r{s_0}t+\mass r{s_1}t}\\
r((s_0+s_1)t)\ar[dr]_{r\rdis{s_0}{s_1}t}&&r(s_0t)+r(s_1t)\\
&r(s_0t+s_1t)\ar[ur]_{\ldis r{s_0t}{s_1t}}%
}
\]

\[
\xymatrix@C=.5em@L=2ex
{
&(r_0s+r_1s)t\ar[rr]^{\rdis{r_0s}{r_1s}t}&&(r_0s)t+(r_1s)t\ar[dr]^{\mass{r_0}st+\mass{r_1}st}\\
((r_0+r_1)s)t\ar[ur]^{\rdis{r_0}{r_1}st}\ar[drr]_{\mass{r_0+r_1}st}&&&&r_0(st)+r_1(st)\\
&&(r_0+r_1)(st)\ar[urr]_{\rdis{r_0}{r_1}st}%
}%
\]

\[
\xymatrix@L=1ex{%
1(r_0+r_1)%
\ar[rr]^{\ldis1{r_0}{r_1}}%
\ar[dr]_{\lambda.(r_0+r_1)}%
&&1r_0+1r_1%
\ar[dl]^{\ \lambda.(r_0)+\lambda.(r_1)}%
\\
&r_0+r_1%
\\
(r_0+r_1)1%
\ar[rr]_{\rdis{r_0}{r_1}1}%
\ar[ur]^{\rho.(r_0+r_1)}%
&&r_01+r_11%
\ar[ul]_{\ \rho.(r_0)+\rho.(r_1)}%
}%
\]

\

\[
\xymatrix@C=.5em
{
&r(s_{00}+s_{01})+r(s_{10}+s_{11})%
\ar[dr]^{\hskip3em\ldis r{s_{00}}{s_{01}}+\ldis r{s_{10}}{s_{11}}}%
\\
r((s_{00}+s_{01})+(s_{10}+s_{11}))%
\ar[ur]^{\ldis r{s_{00}+s_{01}}{s_{10}+s_{11}}\hskip3em}%
\ar[d]_{r\coma{s_{00}}{s_{01}}{s_{10}}{s_{11}}}%
&&(rs_{00}+rs_{01})+(rs_{10}+rs_{11})%
\ar[d]^{\coma{rs_{00}}{rs_{01}}{rs_{10}}{rs_{11}}}%
\\
r((s_{00}+s_{10})+(s_{01}+s_{11}))%
\ar[dr]_{\ldis r{s_{00}+s_{10}}{s_{01}+s_{11}}\hskip3em}%
&&(rs_{00}+rs_{10})+(rs_{01}+rs_{11})%
\\
&r(s_{00}+s_{10})+r(s_{01}+s_{11})%
\ar[ur]_{\hskip3em\ldis r{s_{00}}{s_{10}}+\ldis r{s_{01}}{s_{11}}}%
}%
\]

\[
\xymatrix@C=.5em
{
&&(r_0+r_1)(s_0+s_1)%
\ar[drr]^{\hskip3em\ldis{r_0+r_1}{s_0}{s_1}}%
\ar[dll]_{\rdis{r_0}{r_1}{s_0+s_1}\hskip3em}%
\\
r_0(s_0+s_1)+r_1(s_0+s_1)%
\ar[d]_{\ldis{r_0}{s_0}{s_1}+\ldis{r_1}{s_0}{s_1}}%
&&&&(r_0+r_1)s_0+(r_0+r_1)s_1%
\ar[d]^{\rdis{r_0}{r_1}{s_0}+\rdis{r_0}{r_1}{s_1}}%
\\
(r_0s_0+r_0s_1)+(r_1s_0+r_1s_1)%
\ar[rrrr]^{\coma{r_0s_0}{r_0s_1}{r_1s_0}{r_1s_1}}%
&&&&(r_0s_0+r_1s_0)+(r_0s_1+r_1s_1)%
}%
\]

\[
\xymatrix@C=.5em
{
&(r_{00}+r_{01})s+(r_{10}+r_{11})s%
\ar[dr]^{\hskip3em\rdis{r_{00}}{r_{01}}s+\rdis{r_{10}}{r_{11}}s}%
\\
((r_{00}+r_{01})+(r_{10}+r_{11}))s%
\ar[ur]^{\rdis{r_{00}+r_{01}}{r_{10}+r_{11}}s\hskip3em}%
\ar[d]_{\coma{r_{00}}{r_{01}}{r_{10}}{r_{11}}s}%
&&(r_{00}s+r_{01}s)+(r_{10}s+r_{11}s)%
\ar[d]^{\coma{r_{00}s}{r_{01}s}{r_{10}s}{r_{11}s}}%
\\
((r_{00}+r_{10})+(r_{01}+r_{11}))s%
\ar[dr]_{\rdis{r_{00}+r_{10}}{r_{01}+r_{11}}s\hskip3em}%
&&(r_{00}s+r_{10}s)+(r_{01}s+r_{11}s)%
\\
&(r_{00}+r_{10})s+(r_{01}+r_{11})s%
\ar[ur]_{\hskip3em\rdis{r_{00}}{r_{10}}s+\rdis{r_{01}}{r_{11}}s}%
}%
\]
\end{Definition}

Morphisms of categorical rings are defined as follows:

\begin{Definition}
A \emph{2-homomorphism} $\fb:\sR\to\sR'$ is a quadruple
$(f,\fp,\fm,f_1)$ where $f$ is a functor from $\sR$ to $\sR'$,
$\fp$, $\fm$ are natural morphisms of the form
\begin{align*}
\fp(r_0,r_1)&:f(r_0)+f(r_1)\to f(r_0+r_1),\\
\fm(r,s)&:f(r)f(s)\to f(rs)
\end{align*}
and $f_1:f(1_\sR)\to1_{\sR'}$ is a morphism such that $(f,\fp,f_0)$
and $(f,\fm,f_1)$ are monoidal functor structures with respect to
the monoidal structures corresponding to $+$ and $\cdot$
respectively, with $f_0$ as in \eqref{zero}, and moreover the
diagrams
\[
\xymatrix@!C=7em@L=2ex{%
&f(r)f(s_0)+f(r)f(s_1)%
\ar[r]^{\fm(r,s_0)+\fm(r,s_1)}%
&f(rs_0)+f(rs_1)%
\ar[dr]^{\fp(rs_0,rs_1)}%
\\
f(r)(f(s_0)+f(s_1))%
\ar[ur]^{\ldis{f(r)}{f(s_0)}{f(s_1)}\qquad}%
\ar[dr]_{f(r)\fp(s_0,s_1)}%
&&&f(rs_0+rs_1)%
\\
&f(r)f(s_0+s_1)%
\ar[r]_{\fm(r,s_0+s_1)}%
&f(r(s_0+s_1))%
\ar[ur]_{f\left(\ldis r{s_0}{s_1}\right)}%
}
\]
and
\[
\xymatrix@!C=7em@L=2ex{%
&f(r_0)f(s)+f(r_1)f(s)%
\ar[r]^{\fm(r_0,s)+\fm(r_1,s)}%
&f(r_0s)+f(r_1s)%
\ar[dr]^{\fp(r_0s,r_1s)}%
\\
(f(r_0)+f(r_1))f(s)%
\ar[ur]^{\rdis{f(r_0)}{f(r_1)}{f(s)}\qquad}%
\ar[dr]_{\fp(r_0,r_1)f(s)}%
&&&f(r_0s+r_1s)%
\\
&f(r_0+r_1)f(s)%
\ar[r]_{\fm(r_0+r_1,s)}%
&f((r_0+r_1)s)%
\ar[ur]_{f\left(\rdis{r_0}{r_0}s\right)}%
}
\]
commute for all possible objects involved.
\end{Definition}

We will need the following fact in what follows:

\begin{Proposition}\label{holem}
In any categorical ring $\sR$ one has
\[
\xymatrix@C=4em{%
x%
\ar@/^1pc/[r]^{\alpha}%
\ar@{}[r]|{\beta}%
\ar@/_1pc/[r]_{\alpha'}%
&y%
}%
\qquad\then\qquad%
\xymatrix@C=4em{%
rx%
\ar@/^1pc/[r]^{r\alpha}%
\ar@{}[r]|{r\beta}%
\ar@/_1pc/[r]_{r\alpha'}%
&ry,%
&xr%
\ar@/^1pc/[r]^{\alpha r}%
\ar@{}[r]|{\beta r}%
\ar@/_1pc/[r]_{\alpha'r}%
&yr.%
}%
\]
\end{Proposition}

\begin{proof}
It follows from the definition of categorical ring that for any
$r\in\sR$ the morphisms $\ldis r--$, resp. $\rdis--r$, constitute a
structure of a monoidal functor on the endofunctor $r\cdot$, resp.
$\cdot r:\sR\to\sR$ of $\sR$ (with respect to the additive monoidal
structure). The proposition is thus particular case of
\eqref{holef}.
\end{proof}

\section{Third Mac Lane cohomology group}

We refer to \cites{louvain,ill} for the original construction of
Mac~Lane cohomology. Here we will only give the definition of the
third cohomology group as this is all that we need. To make
expressions shorter, we will need the cross-effect notation. Recall
that for a map $f:A\to B$ between abelian groups, its first
cross-effect is a map
\[
\cref f--:A\x A\to B
\]
is given by
\[
\cref fxy=f(x)+f(y)-f(x+y).
\]

\begin{Definition}\label{h3}
For a ring $R$ and an $R$-bimodule $B$, the group $C^3(R;B)$ of
Mac~Lane 3-cochains of $R$ with coefficients in $B$ consists of
quadruples $(\phm,\ph_{\cdot+},\ph_{+\cdot},\php)$, of maps
\begin{align*}
\phm,\ph_{\cdot+},\ph_{+\cdot}:R^3\to B%
\intertext{and}%
\php:R^4\to B
\end{align*}
which are \emph{normalized} in the sense that $\phm$, $\ph_{\cdot+}$
and $\ph_{+\cdot}$ take zero values if one of their arguments is
zero, and moreover $\php$ satisfies
\[
\php\left(\smat{r_0&r_1\\0&0}\right)=\php\left(\smat{0&0\\r_0&r_1}\right)=\php\left(\smat{r_0&0\\r_1&0}\right)=\php\left(\smat{0&r_0\\0&r_1}\right)
=\php\left(\smat{r_0&0\\0&r_1}\right)=0
\]
for all $r_0,r_1\in R$. The group structure is given by valuewise
addition of functions.

The subgroup $Z^3(R;B)\subseteq C^3(R;B)$ of 3-cocycles is singled
out by the following equations:
\begin{align*}
r\phm(s,t,u)-\phm(rs,t,u)\quad\\
+\phm(r,st,u)-\phm(r,s,tu)+\phm(r,s,t)u=&0,\\
r\ph_{\cdot+}(s,t_0,t_1)-\ph_{\cdot+}(rs,t_0,t_1)+\ph_{\cdot+}(r,st_0,st_1)=&\cref{\phm(r,s,-)}{t_0}{t_1},\\
\ph_{\cdot+}(r,s_0t,s_1t)-\ph_{\cdot+}(r,s_0,s_1)t=&\ph_{+\cdot}(rs_0,rs_1,t)-r\ph_{+\cdot}(s_0,s_1,t),\\
\ph_{+\cdot}(r_0s,r_1s,t)-\ph_{+\cdot}(r_0,r_1,st)+\ph_{+\cdot}(r_0,r_1,s)t=&-\cref{\phm(-,s,t)}{r_0}{r_1},\\
\php\left(\smat{rs_{00}&rs_{01}\\rs_{10}&rs_{11}}\right)-r\php\left(\smat{s_{00}&s_{01}\\s_{10}&s_{11}}\right)
=&\cref{\ph_{\cdot+}(r,-,-)}{(s_{00},s_{10})}{(s_{01},s_{11})}\\
&-\cref{\ph_{\cdot+}(r,-,-)}{(s_{00},s_{01})}{(s_{10},s_{11})},\\
\php\left(\smat{r_0s_0&r_0s_1\\r_1s_0&r_1s_1}\right)=&\cref{\ph_{\cdot+}(-,s_0,s_1)}{r_0}{r_1}-\cref{\ph_{+\cdot}(r_0,r_1,-)}{s_0}{s_1},\\
\php\left(\smat{r_{00}s&r_{01}s\\r_{10}s&r_{11}s}\right)-\php\left(\smat{r_{00}&r_{01}\\r_{10}&r_{11}}\right)s
=&\cref{\ph_{+\cdot}(-,-,s)}{(r_{00},r_{10})}{(r_{01},r_{11})}\\
&-\cref{\ph_{+\cdot}(-,-,s)}{(r_{00},r_{01})}{(r_{10},r_{11})},\\
\cref\php{\left(\smat{r_{000}&r_{001}\\r_{010}&r_{011}}\right)}{\left(\smat{r_{100}&r_{101}\\r_{110}&r_{111}}\right)}\quad\\
-\cref\php{\left(\smat{r_{000}&r_{001}\\r_{100}&r_{101}}\right)}{\left(\smat{r_{010}&r_{011}\\r_{110}&r_{111}}\right)}\quad\\
+\cref\php{\left(\smat{r_{000}&r_{010}\\r_{100}&r_{110}}\right)}{\left(\smat{r_{001}&r_{011}\\r_{101}&r_{111}}\right)}=&0.
\end{align*}

The subgroup $B^3(R;B)\subseteq Z^3(R;B)$ of 3-coboundaries consists
of those quadruples $(\phm,\ph_{\cdot+},\ph_{+\cdot},\php)$ for
which there exist maps $\gammam,\gammap:R^2\to B$ such that
\begin{align*}
\phm(r,s,t)=&r\gammam(s,t)-\gammam(rs,t)+\gammam(r,st)-\gammam(r,s)t,\\
\ph_{\cdot+}(r,s_0,s_1)=&r\gammap(s_0,s_1)-\gammap(rs_0,rs_1)+\cref{\gammam(r,-)}{s_0}{s_1},\\
\ph_{+\cdot}(r_0,r_1,s)=&\gammap(r_0s,r_1s)-\gammap(r_0,r_1)s-\cref{\gammam(-,s)}{r_0}{r_1},\\
\php\left(\smat{r_{00}&r_{01}\\r_{10}&r_{11}}\right)=&\cref{\gammap}{(r_{00},r_{01})}{(r_{10},r_{11})}-\cref{\gammap}{(r_{00},r_{10})}{(r_{01},r_{11})}
\end{align*}
for all $r,...\in R$.

Finally, we define
\[
H^3(R;B):=Z^3(R;B)/B^3(R;B).
\]
\end{Definition}

\section{Characteristic class of a categorical ring}\label{ch}

Suppose given a categorical ring $\sR$ as in \ref{2ring}. Then the
set $R=\pi_0(\sR)$ of isomorphism classes of objects of $\sR$ is a
ring, and the group $B=\pi_1(\sR)$ of automorphisms of the zero
object of $\sR$ is an $\sR$-bimodule. We are going to assign to
$\sR$ a cohomology class
\[
\brk{\sR}\in H^3(\pi_0(\sR);\pi_1(\sR)).
\]
For that, we arbitrarily choose an object $\bar r$ of $\sR$ in each
isomorphism class $r\in\pi_0(\sR)$; moreover we arbitrarily choose
morphisms
\[
\sigmam(r,s):\sbar r\sbar s\to\ol{rs}
\]
and
\[
\sigmap(r_0,r_1):\sbar r_0+\sbar r_1\to\ol{r_0+r_1}.
\]
These morphisms give rise to several not necessarily commutative
diagrams. They define elements of $\pi_1(\sR)$ as in section
\ref{catg}.

In particular, for any $r,s,t\in\pi_0(\sR)$ the diagram
\[
\xymatrix@!C=2em{%
&\sbar r(\sbar s\sbar t)%
\ar[rr]^{\sbar r\sigmam(s,t)}%
&&\sbar r\ol{st}%
\ar[dr]|{\sigmam(r,st)}%
\\
(\sbar r\sbar s)\sbar t%
\ar[ur]^{\mass{\sbar r}{\sbar s}{\sbar t}}%
\ar[drr]|{\sigmam(r,s)\sbar t}%
\ar@{}[rrrr]|{\phm(r,s,t)}%
&&&&\ol{rst}%
\\
&&\ol{rs}\,\sbar t%
\ar[urr]|{\sigmam(rs,t)}%
}
\]
produces an element $\phm(r,s,t)\in\pi_1(\sR)$ measuring deviation
from its commutativity. Similarly for any $r,s_0,s_1\in\pi_0(\sR)$
deviation from commutativity of the diagram
\[
\xymatrix@C=3em{%
&\sbar r\sbar s_0+\sbar r\sbar s_1%
\ar[rr]^-{\sigmam(r,s_0)+\sigmam(r,s_1)}%
&&\ol{rs}_0+\ol{rs}_1%
\ar[dr]|{\sigmap(rs_0,rs_1)}%
\\
\sbar r(\sbar s_0+\sbar s_1)%
\ar[ur]^<<{\ldis{\sbar r}{\sbar s_0}{\sbar s_1}}%
\ar[dr]|{\sbar r\sigmap(s_0,s_1)}%
\ar@{}[rrrr]|{\ph_{\cdot+}(r,s_0,s_1)}%
&&&&\ol{rs_0+rs_1}%
\\
&\sbar r\ol{s_0+s_1}%
\ar[rr]_{\sigmam(r,s_0+s_1)}%
&&\ol{r(s_0+s_1)}%
\ar@{=}[ur]%
}
\]
is measured by an element $\ph_{\cdot+}(r,s_0,s_1)\in\pi_1(\sR)$;
for any $r_0,r_1,s\in\pi_0(\sR)$ the diagram
\[
\xymatrix@C=3em{%
&\sbar r_0\sbar s+\sbar r_1\sbar s%
\ar[rr]^-{\sigmam(r_0,s)+\sigmam(r_1,s)}%
&&\ol{r_0s}+\ol{r_1s}%
\ar[dr]|{\sigmap(r_0s,r_1s)}%
\\
(\sbar r_0+\sbar r_1)\sbar s%
\ar[ur]^<<{\rdis{\sbar r_0}{\sbar r_1}{\sbar s}}%
\ar[dr]|{\sigmap(r_0,r_1)\sbar s}%
\ar@{}[rrrr]|{\ph_{+\cdot}(r_0,r_1,s)}%
&&&&\ol{r_0s+r_1s}%
\\
&\ol{r_0+r_1}\,\sbar s%
\ar[rr]_{\sigmam(r_0+r_1,s)}%
&&\ol{(r_0+r_1)s}%
\ar@{=}[ur]%
}
\]
gives $\ph_{+\cdot}(r_0,r_1,s)\in\pi_1(\sR)$; and for any
$r_{00},r_{01},r_{10},r_{11}\in\pi_0(\sR)$ the diagram
\[
\xymatrix@!C=10em@R=3em{%
&\ol{r_{00}+r_{01}}+\ol{r_{10}+r_{11}}%
\ar[dr]|{\sigmap(r_{00}+r_{01},r_{10}+r_{11})}%
\ar@{}[ddd]|{\php\left(\smat{r_{00}&r_{01}\\r_{10}&r_{11}}\right)}%
\\
(\sbar r_{00}+\sbar r_{01})+(\sbar r_{10}+\sbar r_{11})%
\ar[ur]|{\sigmap(r_{00},r_{01})+\sigmap(r_{10},r_{11})}%
\ar[d]_{\coma{\sbar r_{00}}{\sbar r_{01}}{\sbar r_{10}}{\sbar r_{11}}}%
&&\ol{r_{00}+r_{01}+r_{10}+r_{11}}%
\ar@{=}[d]%
\\
(\sbar r_{00}+\sbar r_{10})+(\sbar r_{01}+\sbar r_{11})%
\ar[dr]|{\sigmap(r_{00},r_{10})+\sigmap(r_{01},r_{11})}%
&&\ol{r_{00}+r_{10}+r_{01}+r_{11}}%
\\
&\ol{r_{00}+r_{10}}+\ol{r_{01}+r_{11}}%
\ar[ur]|{\sigmap(r_{00}+r_{10},r_{01}+r_{11})}%
}
\]
gives
$\php\left(\smat{r_{00}&r_{01}\\r_{10}&r_{11}}\right)\in\pi_1(\sR)$.

Thus the above diagrams give rise to a 3-cochain $\ph$ in the
Mac~Lane complex of $\pi_0(\sR)$ with coefficients in $\pi_1(\sR)$.
Explicitly, it is defined by
\begin{align*}
\phm(r,s,t)&=\sigmam(r,st)\circ\sbar r\sigmam(s,t)\circ\mass{\sbar r}{\sbar s}{\sbar t}-\sigmam(rs,t)\circ\sigmam(r,s)\sbar t,\\
\ph_{\cdot+}(r,s_0,s_1)&=\sigmap(rs_0,rs_1)\circ(\sigmam(r,s_0)+\sigmam(r,s_1))\circ\ldis{\sbar r}{\sbar s_0}{\sbar s_1}-\sigmam(r,s_0+s_1)\circ\sbar r\sigmap(s_0,s_1),\\
\ph_{+\cdot}(r_0,r_1,s)&=\sigmap(r_0s,r_1s)\circ(\sigmam(r_0,s)+\sigmam(r_1,s))\circ\rdis{\sbar r_0}{\sbar r_1}{\sbar s}-\sigmam(r_0+r_1,s)\circ\sigmap(r_0,r_1)\sbar s,\\
\php\left(\smat{r_{00}&r_{01}\\r_{10}&r_{11}}\right)&=\sigmap(r_{00}+r_{10},r_{01}+r_{11})\circ(\sigmap(r_{00},r_{10})+\sigmap(r_{01},r_{11}))\circ\coma{\sbar r_{00}}{\sbar r_{01}}{\sbar r_{10}}{\sbar r_{11}}\\
&-\sigmap(r_{00}+r_{01},r_{10}+r_{11})\circ(\sigmap(r_{00},r_{01})+\sigmap(r_{10},r_{11})).
\end{align*}

We have

\begin{Proposition}
The above cochain is a cocycle.
\end{Proposition}

\begin{proof}
We have to check the eight equalities from \eqref{h3}. These
equalities are deducible from considering eight diagrams below. In
all of these diagrams, all quadrangles commute by naturality, the
inner pentagons are filled by the indicated elements of $\pi_1(\sR)$
using \eqref{holep} and \eqref{holem} as needed, whereas the outer
perimeters commute since each of them coincides with a coherence
diagram from the definition of categorical ring.

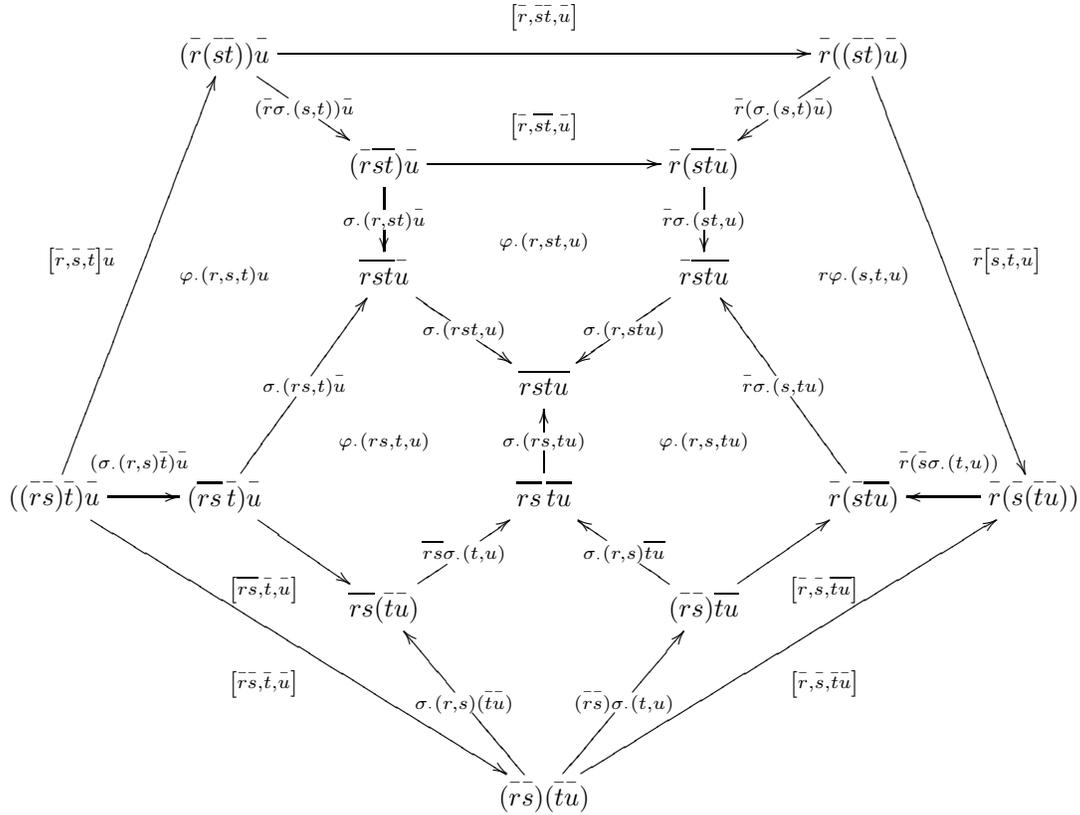
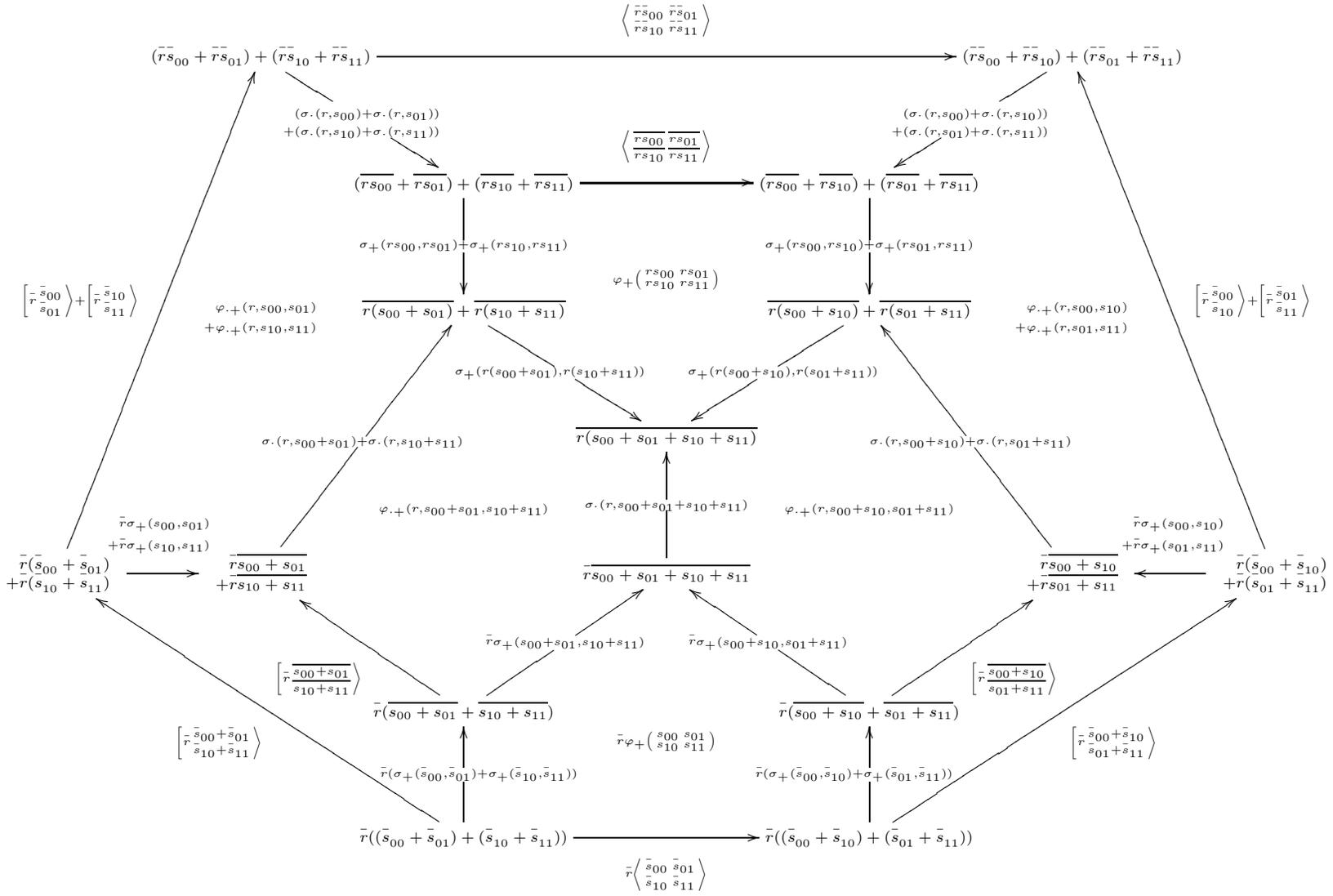
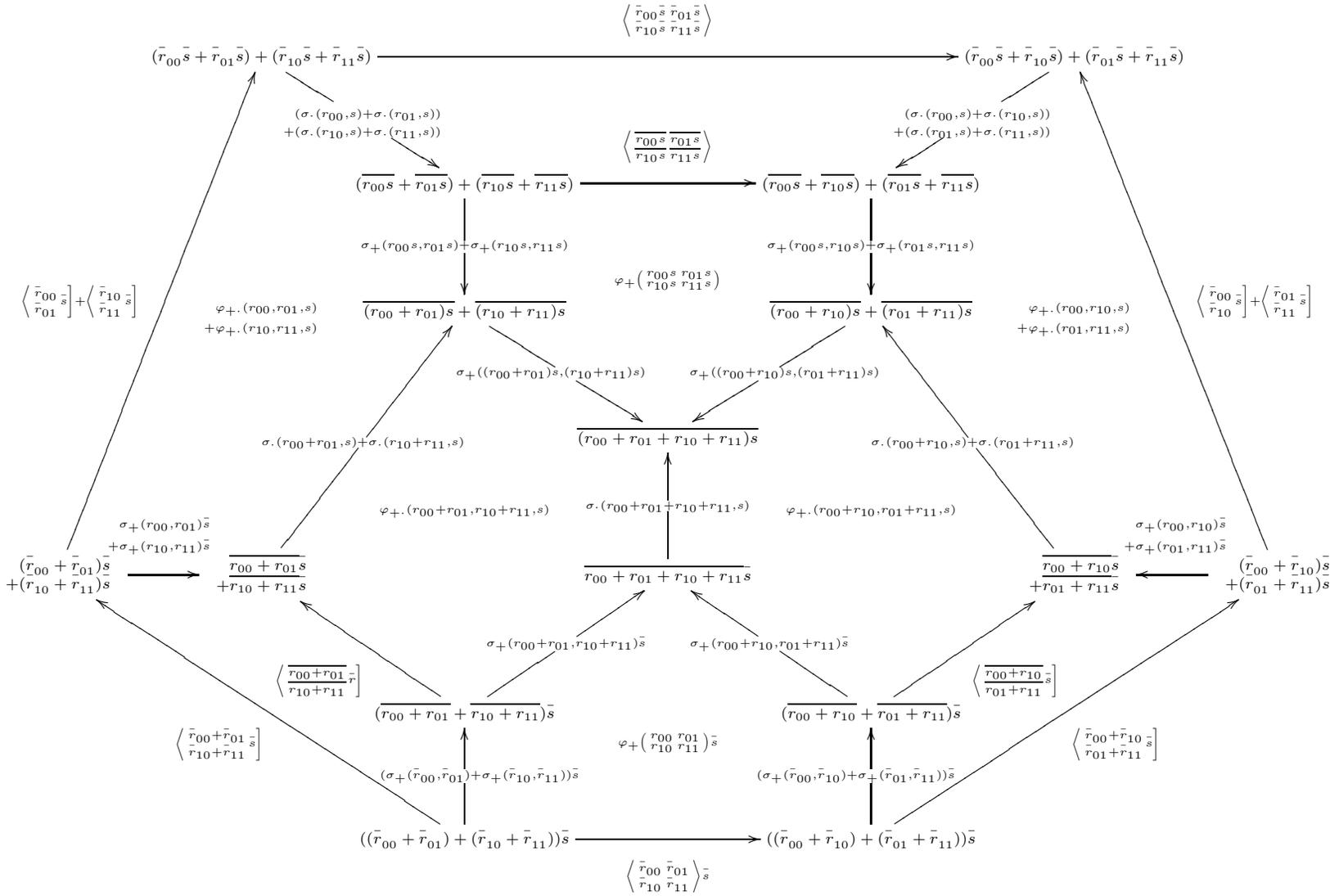
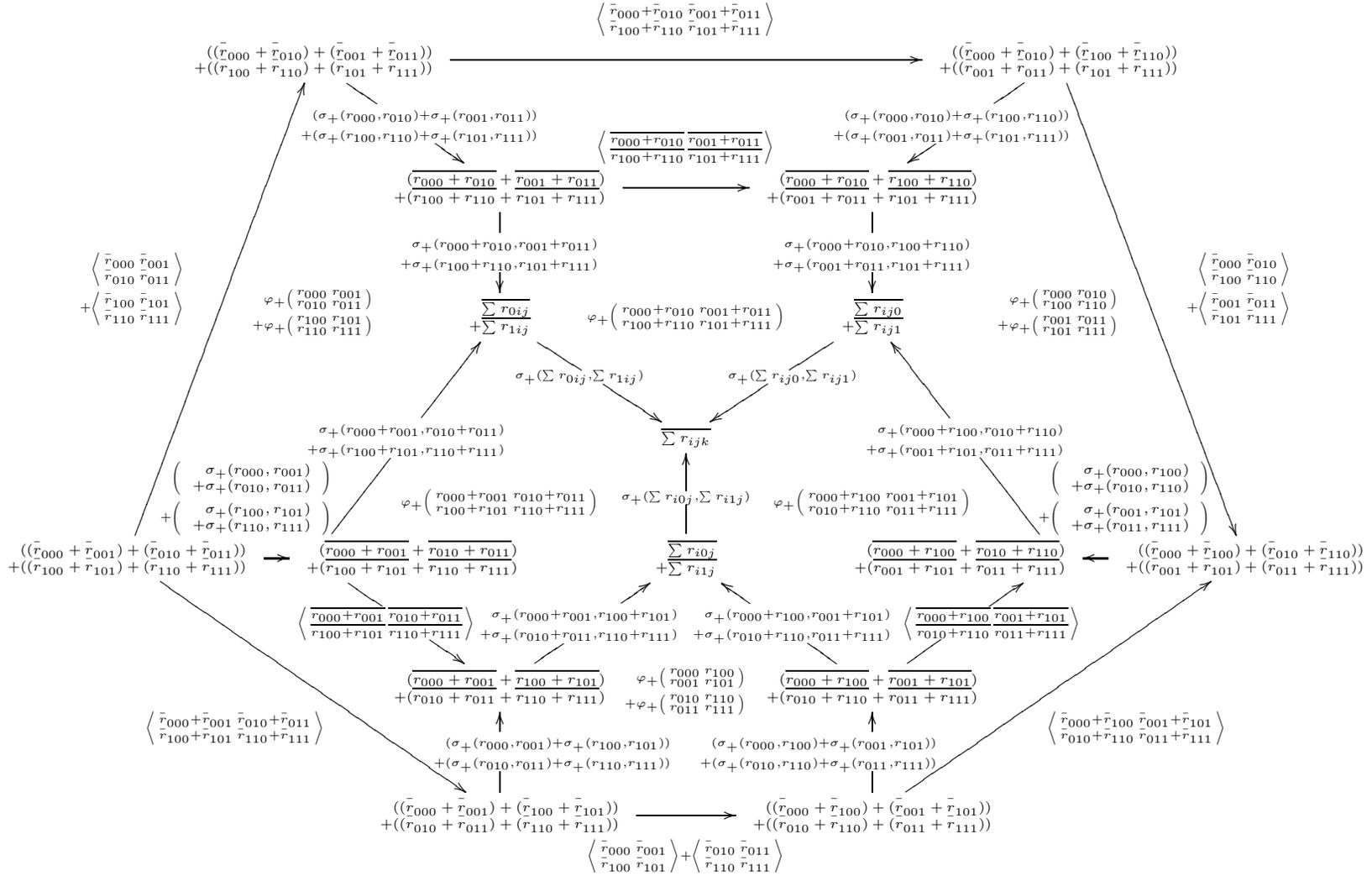
\[
\xymatrix@L=2ex{%
&(\sbar r(\sbar s\sbar t))\sbar u%
\ar@{}[dddd]|{\phm(r,s,t)u}%
\ar[dr]|{(\sbar r\sigmam(s,t))\sbar u}%
\ar[rrrr]^{\mass{\sbar r}{\sbar s\sbar t}{\sbar u}}%
&&&&\sbar r((\sbar s\sbar t)\sbar u)%
\ar@{}[dddd]|{r\phm(s,t,u)}%
\ar[dl]|{\sbar r(\sigmam(s,t)\sbar u)}%
\ar[ddddr]^{\sbar r\mass{\sbar s}{\sbar t}{\sbar u}}%
\\
&&(\sbar r\ol{st})\sbar u%
\ar[d]|{\sigmam(r,st)\sbar u}%
\ar[rr]^{\mass{\sbar r}{\ol{st}}{\sbar u}}%
&&\sbar r(\ol{st}\sbar u)%
\ar[d]|{\sbar r\sigmam(st,u)}%
\\
&&\ol{rst}\sbar u%
\ar@{}[rr]^{\phm(r,st,u)}%
\ar@{}[ddd]|{\phm(rs,t,u)}%
\ar[dr]|{\sigmam(rst,u)}%
&&\sbar r\ol{stu}%
\ar@{}[ddd]|{\phm(r,s,tu)}%
\ar[dl]|{\sigmam(r,stu)}%
\\
&&&\ol{rstu}%
\\
((\sbar r\sbar s)\sbar t)\sbar u%
\ar[r]^{(\sigmam(r,s)\sbar t)\sbar u}%
\ar[dddrrr]_{\mass{\sbar r\sbar s}{\sbar t}{\sbar u}}%
\ar[uuuur]^{\mass{\sbar r}{\sbar s}{\sbar t}\sbar u}%
&(\ol{rs}\,\sbar t)\sbar u%
\ar[dr]_{\mass{\ol{rs}}{\sbar t}{\sbar u}}%
\ar[uur]|{\sigmam(rs,t)\sbar u}%
&&\ol{rs}\,\ol{tu}%
\ar[u]|{\sigmam(rs,tu)}%
&&\sbar r(\sbar s\ol{tu})%
\ar[uul]|{\sbar r\sigmam(s,tu)}%
&\sbar r(\sbar s(\sbar t\sbar u))%
\ar[l]_{\sbar r(\sbar s\sigmam(t,u))}%
\\
&&\ol{rs}(\sbar t\sbar u)%
\ar[ur]|{\ol{rs}\sigmam(t,u)}%
&&(\sbar r\sbar s)\ol{tu}%
\ar[ul]|{\sigmam(r,s)\ol{tu}}%
\ar[ur]_{\mass{\sbar r}{\sbar s}{\ol{tu}}}%
\\
\\
&&&(\sbar r\sbar s)(\sbar t\sbar u)%
\ar[uul]|{\sigmam(r,s)(\sbar t\sbar u)}%
\ar[uur]|{(\sbar r\sbar s)\sigmam(t,u)}%
\ar[uuurrr]_{\mass{\sbar r}{\sbar s}{\sbar t\sbar u}}%
}
\]

\begin{landscape}

\[
\xymatrix@R=3.6em@L=2ex{%
&\sbar r(\sbar s\sbar t_0+\sbar s\sbar t_1)%
\ar@{}[dddd]|{r\ph_{\cdot+}(s,t_0,t_1)}%
\ar[dr]|{\sbar r(\sigmam(s,t_0)+\sigmam(s,t_1))}%
\ar[rrrr]^{\ldis{\sbar r}{\sbar s\sbar t_0}{\sbar s\sbar t_1}}%
&&&&\sbar r(\sbar s\sbar t_0)+\sbar r(\sbar s\sbar t_1)%
\ar@{}[dddd]|{\phm(r,s,t_0)+\phm(r,s,t_1)}%
\ar[dl]|{\sbar r\sigmam(s,t_0)+\sbar r\sigmam(s,t_1)}%
\\
&&\sbar r(\ol{st_0}+\ol{st_1})%
\ar[d]|{\sbar r\sigmap(st_0,st_1)}%
\ar[rr]^{\ldis{\sbar r}{\ol{st_0}}{\ol{st_1}}}%
&&\sbar r\ol{st_0}+\sbar r\ol{st_1}%
\ar[d]|{\sigmam(r,st_0)+\sigmam(r,st_1)}%
\\
&&\sbar r\ol{s(t_0+t_1)}%
\ar@{}[rr]^{\ph_{\cdot+}(r,st_0,st_1)}%
\ar@{}[ddd]|{\phm(r,s,t_0+t_1)}%
\ar[dr]|{\sigmam(r,s(t_0+t_1))}%
&&\ol{rst_0}+\ol{rst_1}%
\ar@{}[ddd]|{\ph_{\cdot+}(rs,t_0,t_1)}%
\ar[dl]|{\sigmap(rst_0,rst_1)}%
\\
&&&\ol{rs(t_0+t_1)}%
\\
\sbar r(\sbar s(\sbar t_0+\sbar t_1))%
\ar[r]^{\sbar r(\sbar s\sigmap(t_0,t_1))}%
\ar[uuuur]^{\sbar r\ldis{\sbar s}{\sbar t_0}{\sbar t_1}}%
&\sbar r(\sbar s\ol{t_0+t_1})%
\ar[uur]|{\sbar r\sigmam(s,t_0+t_1)}%
&&\ol{rs}\,\ol{t_0+t_1}%
\ar[u]|{\sigmam(rs,t_0+t_1)}%
&&\ol{rs}\,\sbar t_0+\ol{rs}\,\sbar t_1%
\ar[uul]|{\sigmam(rs,t_0)+\sigmam(rs,t_1)}%
&(\sbar r\sbar s)\sbar t_0+(\sbar r\sbar s)\sbar t_1%
\ar[l]_{\sigmam(r,s)\sbar t_0+\sigmam(r,s)\sbar t_1}%
\ar[uuuul]_{\mass{\sbar r}{\sbar s}{\sbar t_0}+\mass{\sbar r}{\sbar s}{\sbar t_1}}%
\\
&&(\sbar r\sbar s)\ol{t_0+t_1}%
\ar[ur]|{\sigmam(r,s)\ol{t_0+t_1}}%
\ar[ul]^{\mass{\sbar r}{\sbar s}{\ol{t_0+t_1}}}%
&&\ol{rs}(\sbar t_0+\sbar t_1)%
\ar[ul]|{\ol{rs}\sigmap(t_0,t_1)}%
\ar[ur]_{\ldis{\ol{rs}}{\sbar t_0}{\sbar t_1}}%
\\
\\
&&&(\sbar r\sbar s)(\sbar t_0+\sbar t_1)%
\ar[uuulll]^{\mass{\sbar r}{\sbar s}{\sbar t_0+\sbar t_1}}%
\ar[uul]|{(\sbar r\sbar s)\sigmap(t_0,t_1)}%
\ar[uur]|{\sigmam(r,s)(\sbar t_0+\sbar t_1)}%
\ar[uuurrr]_{\ldis{\sbar r\sbar s}{\sbar t_0}{\sbar t_1}}%
}
\]

\[
\xymatrix@R=4.2em@L=2ex{%
&(\sbar r(\sbar s_0+\sbar s_1))\sbar t%
\ar@{}[dddd]|{\ph_{\cdot+}(r,s_0,s_1)t}%
\ar[dr]|{(\sbar r\sigmap(s_0,s_1))\sbar t}%
\ar[rrrr]^{\mass{\sbar r}{\sbar s_0+\sbar s_1}{\sbar t}}%
\ar[ddddl]_{\ldis{\sbar r}{\sbar s_0}{\sbar s_1}\sbar t}%
&&&&\sbar r((\sbar s_0+\sbar s_1)\sbar t)%
\ar@{}[dddd]|{r\ph_{+\cdot}(s_0,s_1,t)}%
\ar[dl]|{\sbar r(\sigmap(s_0,s_1)\sbar t)}%
\ar[ddddr]^{\sbar r\rdis{\sbar s_0}{\sbar s_1}{\sbar t}}%
\\
&&(\sbar r\ol{s_0+s_1})\sbar t%
\ar[d]|{\sigmam(r,s_0+s_1)\sbar t}%
\ar[rr]^{\mass{\sbar r}{\ol{s_0+s_1}}{\sbar t}}%
&&\sbar r(\ol{s_0+s_1}\,\sbar t)%
\ar[d]|{\sbar r\sigmap(s_0+s_1,t)}%
\\
&&\ol{r(s_0+s_1)}\,\sbar t%
\ar@{}[rr]^{\phm(r,s_0+s_1,t)}%
\ar@{}[ddd]|{\ph_{+\cdot}(rs_0,rs_1,t)}%
\ar[dr]|{\sigmam(r(s_0+s_1),t)}%
&&\sbar r\ol{(s_0+s_1)t}%
\ar@{}[ddd]|{\ph_{\cdot+}(r,s_0t,s_1t)}%
\ar[dl]|{\sigmam(r,(s_0+s_1)t)}%
\\
&&&\ol{r(s_0+s_1)t}%
\\
(\sbar r\sbar s_0+\sbar r\sbar s_1)\sbar t%
\ar[r]^{\hskip2em(\sigmam(r,s_0)+\sigmam(r,s_1))\sbar t}%
\ar[ddrr]_{\rdis{\sbar r\sbar s_0}{\sbar r\sbar s_1}{\sbar t}}%
&(\ol{rs_0}+\ol{rs_1})\sbar t%
\ar[dr]_{\rdis{\ol{rs_0}}{\ol{rs_1}}{\sbar t}}%
\ar[uur]|{\sigmap(rs_0,rs_1)\sbar t}%
&&\ol{rs_0t}+\ol{rs_1t}%
\ar[u]|{\sigmap(rs_0t,rs_1t)}%
&&\sbar r(\ol{s_0t}+\ol{s_1t})%
\ar[dl]^{\ldis{\sbar r}{\ol{s_0t}}{\ol{s_1t}}}%
\ar[uul]|{\sbar r\sigmap(s_0t,s_1t)}%
&\sbar r(\sbar s_0\sbar t+\sbar s_1\sbar t)%
\ar[l]_{\sbar r(\sigmam(s_0,t)+\sigmam(s_1,t))\hskip2em}%
\ar[ddll]^{\ldis{\sbar r}{\sbar s_0\sbar t}{\sbar s_1\sbar t}}%
\\
&&\ol{rs_0}\,\sbar t+\ol{rs_1}\,\sbar t%
\ar[ur]|{\sigmam(rs_0,t)+\sigmam(rs_1,t)}%
\ar@{}[rr]_{\phm(r,s_0,t)+\phm(r,s_1,t)}%
&&\sbar r\ol{s_0t}+\sbar r\ol{s_1t}%
\ar[ul]|{\sigmam(r,s_0t)+\sigmam(r,s_1t)}%
\\
&&(\sbar r\sbar s_0)\sbar t+(\sbar r\sbar s_1)\sbar t%
\ar[u]|{\sigmam(r,s_0)\sbar t+\sigmam(r,s_1)\sbar t}%
\ar[rr]_{\mass{\sbar r}{\sbar s_0}{\sbar t}+\mass{\sbar r}{\sbar s_1}{\sbar t}}%
&&\sbar r(\sbar s_0\sbar t)+\sbar r(\sbar s_1\sbar t)%
\ar[u]|{\sbar r\sigmam(s_0,t)+\sbar r\sigmam(s_1,t)}%
}%
\]

\[
\xymatrix@R=3.6em@L=2ex{%
&(\sbar r_0\sbar s+\sbar r_1\sbar s)\sbar t%
\ar@{}[dddd]|{\ph_{+\cdot}(r_0,r_1,s)t}%
\ar[dr]|{(\sigmam(r_0,s)+\sigmam(r_1,s))\sbar t}%
\ar[rrrr]^{\rdis{\sbar r_0\sbar s}{\sbar r_1\sbar s}{\sbar t}}%
&&&&(\sbar r_0\sbar s)\sbar t+(\sbar r_1\sbar s)\sbar t%
\ar@{}[dddd]|{\phm(r_0,s,t)+\phm(r_1,s,t)}%
\ar[dl]|{\sigmam(r_0,s)\sbar t+\sigmam(r_1,s)\sbar t}%
\ar[ddddr]^{\mass{\sbar r_0}{\sbar s}{\sbar t}+\mass{\sbar
r_1}{\sbar s}{\sbar t}}
\\
&&(\ol{r_0s}+\ol{r_1s})\sbar t%
\ar[d]|{\sigmap(r_0s,r_1s)\sbar t}%
\ar[rr]^{\rdis{\ol{r_0s}}{\ol{r_1s}}{\sbar t}}%
&&\ol{r_0s}\,\sbar t+\ol{r_1s}\,\sbar t%
\ar[d]|{\sigmam(r_0s,t)+\sigmam(r_1s,t)}%
\\
&&\ol{(r_0+r_1)s}\,\sbar t%
\ar@{}[rr]^{\ph_{+\cdot}(r_0s,r_1s,t)}%
\ar@{}[ddd]|{\phm(r_0+r_1,s,t)}%
\ar[dr]|{\sigmam((r_0+r_1)s,t)}%
&&\ol{r_0st}+\ol{r_1st}%
\ar@{}[ddd]|{\ph_{+\cdot}(r_0,r_1,st)}%
\ar[dl]|{\sigmap(r_0st,r_1st)}%
\\
&&&\ol{(r_0+r_1)st}%
\\
((\sbar r_0+\sbar r_1)\sbar s)\sbar t%
\ar[r]^{(\sigmap(r_0,r_1)\sbar s)\sbar t}%
\ar[uuuur]^{\rdis{\sbar r_0}{\sbar r_1}{\sbar s}\sbar t}%
\ar[dddrrr]_{\mass{\sbar r_0+\sbar r_1}{\sbar s}{\sbar t}}%
&(\ol{r_0+r_1}\sbar s)\sbar t%
\ar[uur]|{\sigmam(r_0+r_1,s)\sbar t}%
\ar[dr]_{\mass{\ol{r_0+r_1}}{\sbar s}{\sbar t}}
&&\ol{r_0+r_1}\,\ol{st}%
\ar[u]|{\sigmam(r_0+r_1,st)}%
&&\sbar r_0\ol{st}+\sbar r_1\ol{st}%
\ar[uul]|{\sigmam(r_0,st)+\sigmam(r_1,st)}%
&\sbar r_0(\sbar s\sbar t)+\sbar r_1(\sbar s\sbar t)%
\ar[l]_{\sbar r_0\sigmam(s,t)+\sbar r_1\sigmam(s,t)}%
\\
&&\ol{r_0+r_1}(\sbar s\sbar t)%
\ar[ur]|{\ol{r_0+r_1}\sigmam(s,t)}%
&&(\sbar r_0+\sbar r_1)\ol{st}%
\ar[ul]|{\sigmap(r_0,r_1)\ol{st}}%
\ar[ur]_{\rdis{\sbar r_0}{\sbar r_1}{\ol{st}}}%
\\
\\
&&&(\sbar r_0+\sbar r_1)(\sbar s\sbar t)%
\ar[uul]|{\sigmap(r_0,r_1)(\sbar s\sbar t)}%
\ar[uur]|{(\sbar r_0+\sbar r_1)\sigmam(s,t)}%
\ar[uuurrr]_{\rdis{\sbar r_0}{\sbar r_1}{\sbar s\sbar t}}%
}
\]

{\scriptsize%
\xymatrix@R=6.28em@!C=10em@L=3ex{%
&(\sbar r\sbar s_{00}+\sbar r\sbar s_{01})+(\sbar r\sbar s_{10}+\sbar r\sbar s_{11})%
\ar@{}[dddd]|{\rarr{_{\ph_{\cdot+}(r,s_{00},s_{01})}\\_{+\ph_{\cdot+}(r,s_{10},s_{11})}}}%
\ar[dr]|{\rarr{_{(\sigmam(r,s_{00})+\sigmam(r,s_{01}))}\\_{+(\sigmam(r,s_{10})+\sigmam(r,s_{11}))}}}%
\ar[rrrr]^{\coma{\sbar r\sbar s_{00}}{\sbar r\sbar s_{01}}{\sbar r\sbar s_{10}}{\sbar r\sbar s_{11}}}%
&&&&(\sbar r\sbar s_{00}+\sbar r\sbar s_{10})+(\sbar r\sbar s_{01}+\sbar r\sbar s_{11})%
\ar@{}[dddd]|{\rarr{_{\ph_{\cdot+}(r,s_{00},s_{10})}\\_{+\ph_{\cdot+}(r,s_{01},s_{11})}}}%
\ar[dl]|{\rarr{_{(\sigmam(r,s_{00})+\sigmam(r,s_{10}))}\\_{+(\sigmam(r,s_{01})+\sigmam(r,s_{11}))}}}%
\\
&&(\ol{rs_{00}}+\ol{rs_{01}})+(\ol{rs_{10}}+\ol{rs_{11}})%
\ar[d]|{\sigmap(rs_{00},rs_{01})+\sigmap(rs_{10},rs_{11})}%
\ar[rr]^{\coma{\ol{rs_{00}}}{\ol{rs_{01}}}{\ol{rs_{10}}}{\ol{rs_{11}}}}%
&&(\ol{rs_{00}}+\ol{rs_{10}})+(\ol{rs_{01}}+\ol{rs_{11}})%
\ar[d]|{\sigmap(rs_{00},rs_{10})+\sigmap(rs_{01},rs_{11})}%
\\
&&\ol{r(s_{00}+s_{01})}+\ol{r(s_{10}+s_{11})}%
\ar@{}[rr]^{\php\left(\smat{rs_{00}&rs_{01}\\rs_{10}&rs_{11}}\right)}%
\ar@{}[ddd]|{\ph_{\cdot+}(r,s_{00}+s_{01},s_{10}+s_{11})}%
\ar[dr]|{\sigmap(r(s_{00}+s_{01}),r(s_{10}+s_{11}))\qquad}%
&&\ol{r(s_{00}+s_{10})}+\ol{r(s_{01}+s_{11})}%
\ar@{}[ddd]|{\ph_{\cdot+}(r,s_{00}+s_{10},s_{01}+s_{11})}%
\ar[dl]|{\qquad\sigmap(r(s_{00}+s_{10}),r(s_{01}+s_{11}))}%
\\
&&&\ol{r(s_{00}+s_{01}+s_{10}+s_{11})}%
\\
\rarr{\sbar r(\sbar s_{00}+\sbar s_{01})\\+\sbar r(\sbar s_{10}+\sbar s_{11})}%
\ar[r]^{\rarr{_{\sbar r\sigmap(s_{00},s_{01})}\\_{+\sbar r\sigmap(s_{10},s_{11})}}}%
\ar[uuuur]^{\ldis{\sbar r}{\sbar s_{00}}{\sbar s_{01}}+\ldis{\sbar r}{\sbar s_{10}}{\sbar s_{11}}}%
&\rarr{\sbar r\ol{s_{00}+s_{01}}\\+\sbar r\ol{s_{10}+s_{11}}}%
\ar[uur]|{\sigmam(r,s_{00}+s_{01})+\sigmam(r,s_{10}+s_{11})}%
&&\sbar r\ol{s_{00}+s_{01}+s_{10}+s_{11}}%
\ar[u]|{\sigmam(r,s_{00}+s_{01}+s_{10}+s_{11})}%
&&\rarr{\sbar r\ol{s_{00}+s_{10}}\\+\sbar r\ol{s_{01}+s_{11}}}%
\ar[uul]|{\sigmam(r,s_{00}+s_{10})+\sigmam(r,s_{01}+s_{11})}%
&\rarr{\sbar r(\sbar s_{00}+\sbar s_{10})\\+\sbar r(\sbar s_{01}+\sbar s_{11})}%
\ar[l]_{\rarr{_{\sbar r\sigmap(s_{00},s_{10})}\\_{+\sbar r\sigmap(s_{01},s_{11})}}}%
\ar[uuuul]_{\ldis{\sbar r}{\sbar s_{00}}{\sbar s_{10}}+\ldis{\sbar r}{\sbar s_{01}}{\sbar s_{11}}}%
\\
&&\sbar r(\ol{s_{00}+s_{01}}+\ol{s_{10}+s_{11}})%
\ar[ul]^{\ldis{\sbar r}{\ol{s_{00}+s_{01}}}{\ol{s_{10}+s_{11}}}}%
\ar[ur]|{\sbar r\sigmap(s_{00}+s_{01},s_{10}+s_{11})}%
\ar@{}[rr]_{\sbar r\php\left(\smat{s_{00}&s_{01}\\s_{10}&s_{11}}\right)}%
&&\sbar r(\ol{s_{00}+s_{10}}+\ol{s_{01}+s_{11}})%
\ar[ul]|{\sbar r\sigmap(s_{00}+s_{10},s_{01}+s_{11})}%
\ar[ur]_{\ldis{\sbar r}{\ol{s_{00}+s_{10}}}{\ol{s_{01}+s_{11}}}}%
\\
&&\sbar r((\sbar s_{00}+\sbar s_{01})+(\sbar s_{10}+\sbar s_{11}))%
\ar[u]|{\qquad\sbar r(\sigmap(\sbar s_{00},\sbar s_{01})+\sigmap(\sbar s_{10},\sbar s_{11}))}%
\ar[rr]_{\sbar r\coma{\sbar s_{00}}{\sbar s_{01}}{\sbar s_{10}}{\sbar s_{11}}}%
\ar[uull]^{\ldis{\sbar r}{\sbar s_{00}+\sbar s_{01}}{\sbar s_{10}+\sbar s_{11}}}%
&&\sbar r((\sbar s_{00}+\sbar s_{10})+(\sbar s_{01}+\sbar s_{11}))%
\ar[u]|{\sbar r(\sigmap(\sbar s_{00},\sbar s_{10})+\sigmap(\sbar s_{01},\sbar s_{11}))\qquad}%
\ar[uurr]_{\ldis{\sbar r}{\sbar s_{00}+\sbar s_{10}}{\sbar s_{01}+\sbar s_{11}}}%
}%

\xymatrix@R=5.79em@L=3ex@!C=10em{%
&(\sbar r_0\sbar s_0+\sbar r_0\sbar s_1)+(\sbar r_1\sbar s_0+\sbar r_1\sbar s_1)%
\ar@{}[dddd]|{\rarr{_{\ph_{\cdot+}(r_0,s_0,s_1)}\\_{+\ph_{\cdot+}(r_1,s_0,s_1)}}}%
\ar[dr]|{\rarr{_{(\sigmam(r_0,s_0)+\sigmam(r_0,s_1))}\\_{+(\sigmam(r_1,s_0)+\sigmam(r_1,s_1))}}}%
\ar[rrrr]^{\coma{\sbar r_0\sbar s_0}{\sbar r_0\sbar s_1}{\sbar r_1\sbar s_0}{\sbar r_1\sbar s_1}}%
&&&&(\sbar r_0\sbar s_0+\sbar r_1\sbar s_0)+(\sbar r_0\sbar s_1+\sbar r_1\sbar s_1)%
\ar@{}[dddd]|{\rarr{_{\ph_{+\cdot}(r_0,r_1,s_0)}\\_{+\ph_{+\cdot}(r_0,r_1,s_1)}}}%
\ar[dl]|{\rarr{_{(\sigmam(r_0,s_0)+\sigmam(r_1,s_0))}\\_{+(\sigmam(r_0,s_1)+\sigmam(r_1,s_1))}}}%
\\
&&(\ol{r_0s_0}+\ol{r_0s_1})+(\ol{r_1s_0}+\ol{r_1s_1})%
\ar[d]|{\sigmap(r_0s_0,r_0s_1)+\sigmap(r_1s_0,r_1s_1)}%
\ar[rr]^{\coma{\ol{r_0s_0}}{\ol{r_0s_1}}{\ol{r_1s_0}}{\ol{r_1s_1}}}%
&&(\ol{r_0s_0}+\ol{r_1s_0})+(\ol{r_0s_1}+\ol{r_1s_1})%
\ar[d]|{\sigmap(r_0s_0,r_1s_0)+\sigmap(r_0s_1,r_1s_1)}%
\\
&&\ol{r_0(s_0+s_1)}+\ol{r_1(s_0+s_1)}%
\ar@{}[rr]^{\php\left(\smat{r_0s_0&r_0s_1\\r_1s_0&r_1s_1}\right)}%
\ar@{}[ddd]|{\ph_{+\cdot}(r_0,r_1,s_0+s_1)}%
\ar[dr]|{\sigmap(r_0(s_0+s_1),r_1(s_0+s_1))}%
&&\ol{(r_0+r_1)s_0}+\ol{(r_0+r_1)s_1}%
\ar@{}[ddd]|{\ph_{\cdot+}(r_0+r_1,s_0,s_1)}%
\ar[dl]|{\sigmap((r_0+r_1)s_0,(r_0+r_1)s_1)}%
\\
&&&\ol{(r_0+r_1)(s_0+s_1)}%
\\
\rarr{\sbar r_0(\sbar s_0+\sbar s_1)\\+\sbar r_1(\sbar s_0+\sbar s_1)}%
\ar[r]^{\rarr{_{\sbar r_0\sigmap(s_0,s_1)}\\_{+\sbar r_1\sigmap(s_0,s_1)}}}%
\ar[uuuur]^{\ldis{\sbar r_0}{\sbar s_0}{\sbar s_1}+\ldis{\sbar r_1}{\sbar s_0}{\sbar s_1}}%
&\rarr{\sbar r_0\ol{s_0+s_1}\\+\sbar r_1\ol{s_0+s_1}}%
\ar[uur]|{\sigmam(r_0,s_0+s_1)+\sigmam(r_1,s_0+s_1)}%
&&\ol{r_0+r_1}\,\ol{s_0+s_1}%
\ar[u]|{\sigmam(r_0+r_1,s_0+s_1)}%
&&\rarr{\ol{r_0+r_1}\sbar s_0\\+\ol{r_0+r_1}\sbar s_1}%
\ar[uul]|{\sigmam(r_0+r_1,s_0)+\sigmam(r_0+r_1,s_1)}%
&\rarr{(\sbar r_0+\sbar r_1)\sbar s_0\\+(\sbar r_0+\sbar r_1)\sbar s_1}%
\ar[l]_{\rarr{_{\sigmap(r_0,r_1)\sbar s_0}\\_{+\sigmap(r_0,r_1)\sbar s_1}}}%
\ar[uuuul]_{\rdis{\sbar r_0}{\sbar r_1}{\sbar s_0}+\rdis{\sbar r_0}{\sbar r_1}{\sbar s_1}}%
\\
&&(\sbar r_0+\sbar r_1)\ol{s_0+s_1}%
\ar[ul]^{\rdis{\sbar r_0}{\sbar r_1}{\ol{s_0+s_1}}}%
\ar[ur]|{\sigmap(r_0,r_1)\ol{s_0+s_1}}%
&&\ol{r_0+r_1}(\sbar s_0+\sbar s_1)%
\ar[ul]|{\ol{r_0+r_1}\sigmap(s_0,s_1)}%
\ar[ur]_{\ldis{\ol{r_0+r_1}}{\sbar s_0}{\sbar s_1}}%
\\
\\
&&&(\sbar r_0+\sbar r_1)(\sbar s_0+\sbar s_1)%
\ar[uuulll]^{\rdis{\sbar r_0}{\sbar r_1}{\sbar s_0+\sbar s_1}}%
\ar[uul]|{(\sbar r_0+\sbar r_1)\sigmap(s_0,s_1)}%
\ar[uur]|{\sigmap(r_0,r_1)(\sbar s_0+\sbar s_1)}%
\ar[uuurrr]_{\ldis{\sbar r_0+\sbar r_1}{\sbar s_0}{\sbar s_1}}%
}%

\xymatrix@R=6.28em@!C=10em@L=3ex{%
&(\sbar r_{00}\sbar s+\sbar r_{01}\sbar s)+(\sbar r_{10}\sbar s+\sbar r_{11}\sbar s)%
\ar@{}[dddd]|{\rarr{_{\ph_{+\cdot}(r_{00},r_{01},s)}\\_{+\ph_{+\cdot}(r_{10},r_{11},s)}}}%
\ar[dr]|{\rarr{_{(\sigmam(r_{00},s)+\sigmam(r_{01},s))}\\_{+(\sigmam(r_{10},s)+\sigmam(r_{11},s))}}}%
\ar[rrrr]^{\coma{\sbar r_{00}\sbar s}{\sbar r_{01}\sbar s}{\sbar r_{10}\sbar s}{\sbar r_{11}\sbar s}}%
&&&&(\sbar r_{00}\sbar s+\sbar r_{10}\sbar s)+(\sbar r_{01}\sbar s+\sbar r_{11}\sbar s)%
\ar@{}[dddd]|{\rarr{_{\ph_{+\cdot}(r_{00},r_{10},s)}\\_{+\ph_{+\cdot}(r_{01},r_{11},s)}}}%
\ar[dl]|{\rarr{_{(\sigmam(r_{00},s)+\sigmam(r_{10},s))}\\_{+(\sigmam(r_{01},s)+\sigmam(r_{11},s))}}}%
\\
&&(\ol{r_{00}s}+\ol{r_{01}s})+(\ol{r_{10}s}+\ol{r_{11}s})%
\ar[d]|{\sigmap(r_{00}s,r_{01}s)+\sigmap(r_{10}s,r_{11}s)}%
\ar[rr]^{\coma{\ol{r_{00}s}}{\ol{r_{01}s}}{\ol{r_{10}s}}{\ol{r_{11}s}}}%
&&(\ol{r_{00}s}+\ol{r_{10}s})+(\ol{r_{01}s}+\ol{r_{11}s})%
\ar[d]|{\sigmap(r_{00}s,r_{10}s)+\sigmap(r_{01}s,r_{11}s)}%
\\
&&\ol{(r_{00}+r_{01})s}+\ol{(r_{10}+r_{11})s}%
\ar@{}[rr]^{\php\left(\smat{r_{00}s&r_{01}s\\r_{10}s&r_{11}s}\right)}%
\ar@{}[ddd]|{\ph_{+\cdot}(r_{00}+r_{01},r_{10}+r_{11},s)}%
\ar[dr]|{\sigmap((r_{00}+r_{01})s,(r_{10}+r_{11})s)\qquad}%
&&\ol{(r_{00}+r_{10})s}+\ol{(r_{01}+r_{11})s}%
\ar@{}[ddd]|{\ph_{+\cdot}(r_{00}+r_{10},r_{01}+r_{11},s)}%
\ar[dl]|{\qquad\sigmap((r_{00}+r_{10})s,(r_{01}+r_{11})s)}%
\\
&&&\ol{(r_{00}+r_{01}+r_{10}+r_{11})s}%
\\
\rarr{(\sbar r_{00}+\sbar r_{01})\sbar s\\+(\sbar r_{10}+\sbar r_{11})\sbar s}%
\ar[r]^{\rarr{_{\sigmap(r_{00},r_{01})\sbar s}\\_{+\sigmap(r_{10},r_{11})\sbar s}}}%
\ar[uuuur]^{\rdis{\sbar r_{00}}{\sbar r_{01}}{\sbar s}+\rdis{\sbar r_{10}}{\sbar r_{11}}{\sbar s}}%
&\rarr{\ol{r_{00}+r_{01}}\sbar s\\+\ol{r_{10}+r_{11}}\sbar s}%
\ar[uur]|{\sigmam(r_{00}+r_{01},s)+\sigmam(r_{10}+r_{11},s)}%
&&\ol{r_{00}+r_{01}+r_{10}+r_{11}}\sbar s%
\ar[u]|{\sigmam(r_{00}+r_{01}+r_{10}+r_{11},s)}%
&&\rarr{\ol{r_{00}+r_{10}}\sbar s\\+\ol{r_{01}+r_{11}}\sbar s}%
\ar[uul]|{\sigmam(r_{00}+r_{10},s)+\sigmam(r_{01}+r_{11},s)}%
&\rarr{(\sbar r_{00}+\sbar r_{10})\sbar s\\+(\sbar r_{01}+\sbar r_{11})\sbar s}%
\ar[l]_{\rarr{_{\sigmap(r_{00},r_{10})\sbar s}\\_{+\sigmap(r_{01},r_{11})\sbar s}}}%
\ar[uuuul]_{\rdis{\sbar r_{00}}{\sbar r_{10}}{\sbar s}+\rdis{\sbar r_{01}}{\sbar r_{11}}{\sbar s}}%
\\
&&(\ol{r_{00}+r_{01}}+\ol{r_{10}+r_{11}})\sbar s%
\ar[ul]^{\rdis{\ol{r_{00}+r_{01}}}{\ol{r_{10}+r_{11}}}{\sbar r}}%
\ar[ur]|{\sigmap(r_{00}+r_{01},r_{10}+r_{11})\sbar s}%
\ar@{}[rr]_{\php\left(\smat{r_{00}&r_{01}\\r_{10}&r_{11}}\right)\sbar s}%
&&(\ol{r_{00}+r_{10}}+\ol{r_{01}+r_{11}})\sbar s%
\ar[ul]|{\sigmap(r_{00}+r_{10},r_{01}+r_{11})\sbar s}%
\ar[ur]_{\rdis{\ol{r_{00}+r_{10}}}{\ol{r_{01}+r_{11}}}{\sbar s}}%
\\
&&((\sbar r_{00}+\sbar r_{01})+(\sbar r_{10}+\sbar r_{11}))\sbar s%
\ar[u]|{\qquad(\sigmap(\sbar r_{00},\sbar r_{01})+\sigmap(\sbar r_{10},\sbar r_{11}))\sbar s}%
\ar[rr]_{\coma{\sbar r_{00}}{\sbar r_{01}}{\sbar r_{10}}{\sbar r_{11}}\sbar s}%
\ar[uull]^{\rdis{\sbar r_{00}+\sbar r_{01}}{\sbar r_{10}+\sbar r_{11}}{\sbar s}}%
&&((\sbar r_{00}+\sbar r_{10})+(\sbar r_{01}+\sbar r_{11}))\sbar s%
\ar[u]|{(\sigmap(\sbar r_{00},\sbar r_{10})+\sigmap(\sbar r_{01},\sbar r_{11}))\sbar s\qquad}%
\ar[uurr]_{\rdis{\sbar r_{00}+\sbar r_{10}}{\sbar r_{01}+\sbar r_{11}}{\sbar s}}%
}%
}%

{\tiny%
\xymatrix@!C=10em@R=6em@L=4ex{%
&\rarr{((\sbar r_{000}+\sbar r_{010})+(\sbar r_{001}+\sbar r_{011}))\\+((\sbar r_{100}+\sbar r_{110})+(\sbar r_{101}+\sbar r_{111}))}%
\ar@{}[dddd]|{\rarr{_{\php\left(\smat{r_{000}&r_{001}\\r_{010}&r_{011}}\right)}\\_{+\php\left(\smat{r_{100}&r_{101}\\r_{110}&r_{111}}\right)}}}%
\ar[dr]|{\qquad\rarr{_{(\sigmap(r_{000},r_{010})+\sigmap(r_{001},r_{011}))}\\_{+(\sigmap(r_{100},r_{110})+\sigmap(r_{101},r_{111}))}}}%
\ar[rrrr]^{\coma{\sbar r_{000}+\sbar r_{010}}{\sbar r_{001}+\sbar r_{011}}{\sbar r_{100}+\sbar r_{110}}{\sbar r_{101}+\sbar r_{111}}}%
&&&&\rarr{((\sbar r_{000}+\sbar r_{010})+(\sbar r_{100}+\sbar r_{110}))\\+((\sbar r_{001}+\sbar r_{011})+(\sbar r_{101}+\sbar r_{111}))}%
\ar@{}[dddd]|{\rarr{_{\php\left(\smat{r_{000}&r_{010}\\r_{100}&r_{110}}\right)}\\_{+\php\left(\smat{r_{001}&r_{011}\\r_{101}&r_{111}}\right)}}}%
\ar[dl]|{\rarr{_{(\sigmap(r_{000},r_{010})+\sigmap(r_{100},r_{110}))}\\_{+(\sigmap(r_{001},r_{011})+\sigmap(r_{101},r_{111}))}}\qquad}%
\ar[ddddr]^{\rarr{_{\coma{\sbar r_{000}}{\sbar r_{010}}{\sbar r_{100}}{\sbar r_{110}}}\\_{+\coma{\sbar r_{001}}{\sbar r_{011}}{\sbar r_{101}}{\sbar r_{111}}}}}%
\\
&&\rarr{(\ol{r_{000}+r_{010}}+\ol{r_{001}+r_{011}})\\+(\ol{r_{100}+r_{110}}+\ol{r_{101}+r_{111}})}%
\ar[d]|{\rarr{_{\sigmap(r_{000}+r_{010},r_{001}+r_{011})}\\_{+\sigmap(r_{100}+r_{110},r_{101}+r_{111})}}}%
\ar[rr]^{\coma{\ol{r_{000}+r_{010}}}{\ol{r_{001}+r_{011}}}{\ol{r_{100}+r_{110}}}{\ol{r_{101}+r_{111}}}}%
&&\rarr{(\ol{r_{000}+r_{010}}+\ol{r_{100}+r_{110}})\\+(\ol{r_{001}+r_{011}}+\ol{r_{101}+r_{111}})}%
\ar[d]|{\rarr{_{\sigmap(r_{000}+r_{010},r_{100}+r_{110})}\\_{+\sigmap(r_{001}+r_{011},r_{101}+r_{111})}}}%
\\
&&\rarr{\ol{\sum r_{0ij}}\\+\ol{\sum r_{1ij}}}%
\ar@{}[rr]|{\php\left(\smat{r_{000}+r_{010}&r_{001}+r_{011}\\r_{100}+r_{110}&r_{101}+r_{111}}\right)}%
\ar@{}[ddd]|{\php\left(\smat{r_{000}+r_{001}&r_{010}+r_{011}\\r_{100}+r_{101}&r_{110}+r_{111}}\right)}%
\ar[dr]|{\sigmap(\sum r_{0ij},\sum r_{1ij})\qquad}%
&&\rarr{\ol{\sum r_{ij0}}\\+\ol{\sum r_{ij1}}}%
\ar@{}[ddd]|{\php\left(\smat{r_{000}+r_{100}&r_{001}+r_{101}\\r_{010}+r_{110}&r_{011}+r_{111}}\right)}%
\ar[dl]|{\qquad\sigmap(\sum r_{ij0},\sum r_{ij1})}%
\\
&&&\ol{\sum r_{ijk}}%
\\
\rarr{((\sbar r_{000}+\sbar r_{001})+(\sbar r_{010}+\sbar r_{011}))\\+((\sbar r_{100}+\sbar r_{101})+(\sbar r_{110}+\sbar r_{111}))}%
\ar[r]^{\qquad\qquad\rarr{_{\left(\rarr{\sigmap(r_{000},r_{001})\\+\sigmap(r_{010},r_{011})}\right)}\\_{+\left(\rarr{\sigmap(r_{100},r_{101})\\+\sigmap(r_{110},r_{111})}\right)}}}%
\ar[uuuur]^{\rarr{_{\coma{\sbar r_{000}}{\sbar r_{001}}{\sbar r_{010}}{\sbar r_{011}}}\\_{+\coma{\sbar r_{100}}{\sbar r_{101}}{\sbar r_{110}}{\sbar r_{111}}}}}%
\ar[ddrr]_{\coma{\sbar r_{000}+\sbar r_{001}}{\sbar r_{010}+\sbar r_{011}}{\sbar r_{100}+\sbar r_{101}}{\sbar r_{110}+\sbar r_{111}}}%
&**[r]\rarr{(\ol{r_{000}+r_{001}}+\ol{r_{010}+r_{011}})\\+(\ol{r_{100}+r_{101}}+\ol{r_{110}+r_{111}})}%
\ar[uur]|{\rarr{_{\sigmap(r_{000}+r_{001},r_{010}+r_{011})}\\_{+\sigmap(r_{100}+r_{101},r_{110}+r_{111})}}}%
\ar[dr]|{\coma{\ol{r_{000}+r_{001}}}{\ol{r_{010}+r_{011}}}{\ol{r_{100}+r_{101}}}{\ol{r_{110}+r_{111}}}\qquad\quad}%
&&\rarr{\ol{\sum r_{i0j}}\\+\ol{\sum r_{i1j}}}%
\ar[u]|{\sigmap(\sum r_{i0j},\sum r_{i1j})}%
&&**[l]\rarr{(\ol{r_{000}+r_{100}}+\ol{r_{010}+r_{110}})\\+(\ol{r_{001}+r_{101}}+\ol{r_{011}+r_{111}})}%
\ar[uul]|{\rarr{_{\sigmap(r_{000}+r_{100},r_{010}+r_{110})}\\_{+\sigmap(r_{001}+r_{101},r_{011}+r_{111})}}}%
&\rarr{((\sbar r_{000}+\sbar r_{100})+(\sbar r_{010}+\sbar r_{110}))\\+((\sbar r_{001}+\sbar r_{101})+(\sbar r_{011}+\sbar r_{111}))}%
\ar[l]_{\rarr{_{\left(\rarr{\sigmap(r_{000},r_{100})\\+\sigmap(r_{010},r_{110})}\right)}\\_{+\left(\rarr{\sigmap(r_{001},r_{101})\\+\sigmap(r_{011},r_{111})}\right)}}\qquad\qquad}%
\\
&&\rarr{(\ol{r_{000}+r_{001}}+\ol{r_{100}+r_{101}})\\+(\ol{r_{010}+r_{011}}+\ol{r_{110}+r_{111}})}%
\ar[ur]|{\rarr{_{\sigmap(r_{000}+r_{001},r_{100}+r_{101})}\\_{+\sigmap(r_{010}+r_{011},r_{110}+r_{111})}}\qquad}%
\ar@{}[rr]|{\rarr{_{\php\left(\smat{r_{000}&r_{100}\\r_{001}&r_{101}}\right)}\\_{+\php\left(\smat{r_{010}&r_{110}\\r_{011}&r_{111}}\right)}}}%
&&\rarr{(\ol{r_{000}+r_{100}}+\ol{r_{001}+r_{101}})\\+(\ol{r_{010}+r_{110}}+\ol{r_{011}+r_{111}})}%
\ar[ul]|{\qquad\rarr{_{\sigmap(r_{000}+r_{100},r_{001}+r_{101})}\\_{+\sigmap(r_{010}+r_{110},r_{011}+r_{111})}}}%
\ar[ur]|{\qquad\qquad\coma{\ol{r_{000}+r_{100}}}{\ol{r_{001}+r_{101}}}{\ol{r_{010}+r_{110}}}{\ol{r_{011}+r_{111}}}}%
\\
&&\rarr{((\sbar r_{000}+\sbar r_{001})+(\sbar r_{100}+\sbar r_{101}))\\+((\sbar r_{010}+\sbar r_{011})+(\sbar r_{110}+\sbar r_{111}))}%
\ar[u]|{\qquad\qquad\qquad\qquad\rarr{_{(\sigmap(r_{000},r_{001})+\sigmap(r_{100},r_{101}))}\\_{+(\sigmap(r_{010},r_{011})+\sigmap(r_{110},r_{111}))}}}%
\ar[rr]_{\coma{\sbar r_{000}}{\sbar r_{001}}{\sbar r_{100}}{\sbar r_{101}}+\coma{\sbar r_{010}}{\sbar r_{011}}{\sbar r_{110}}{\sbar r_{111}}}%
&&\rarr{((\sbar r_{000}+\sbar r_{100})+(\sbar r_{001}+\sbar r_{101}))\\+((\sbar r_{010}+\sbar r_{110})+(\sbar r_{011}+\sbar r_{111}))}%
\ar[u]|{\rarr{_{(\sigmap(r_{000},r_{100})+\sigmap(r_{001},r_{101}))}\\_{+(\sigmap(r_{010},r_{110})+\sigmap(r_{011},r_{111}))}}\qquad\qquad\qquad\qquad}%
\ar[uurr]_{\coma{\sbar r_{000}+\sbar r_{100}}{\sbar r_{001}+\sbar r_{101}}{\sbar r_{010}+\sbar r_{110}}{\sbar r_{011}+\sbar r_{111}}}%
}%
}%

\end{landscape}

\end{proof}

To show that the above rule indeed determines an assignment of a
cohomology class to each categorical ring, we must also show

\begin{Proposition}\label{propcorr}
A different choice of representative objects $r\mapsto\stil r$ and
morphisms
\begin{align*}
\sigmam'(r,s)&:\stil r\stil s\to\wt{rs},\\
\sigmap'(r_0,r_1)&:\stil r_0+\stil r_1\to\wt{r_0+r_1}
\end{align*}
leads to a cocycle $\ph'$ which is cohomologous to $\ph$.
\end{Proposition}

\begin{proof}
By \eqref{h3}, 3-cocycles $\ph$ and $\ph'$ are cohomologous if and
only if there exist maps $\gammam,\gammap:R\x R\to B$ such that the
following four equalities
\begin{align*}
\ph'_\cdot(r,s,t)=&\phm(r,s,t)+r\gammam(s,t)-\gammam(rs,t)+\gammam(r,st)-\gammam(r,s)t,\\
\ph'_{\cdot+}(r,s_0,s_1)=&\ph_{\cdot+}(r,s_0,s_1)+r\gammap(s_0,s_1)-\gammap(rs_0,rs_1)+\cref{\gammam(r,-)}{s_0}{s_1},\\
\ph'_{+\cdot}(r_0,r_1,s)=&\ph_{+\cdot}(r_0,r_1,s)+\gammap(r_0s,r_1s)-\gammap(r_0,r_1)s-\cref{\gammam(-,s)}{r_0}{r_1},\\
\ph'_+\left(\smat{r_{00}&r_{01}\\r_{10}&r_{11}}\right)=&\php\left(\smat{r_{00}&r_{01}\\r_{10}&r_{11}}\right)+\cref{\gammap}{(r_{00},r_{01})}{(r_{10},r_{11})}-\cref{\gammap}{(r_{00},r_{10})}{(r_{01},r_{11})}
\end{align*}
are satisfied for all possible elements $r,...$ of $R$.

Let us then choose arbitrary morphisms
\[
\shat r:\sbar r\to\stil r
\]
for all $r\in R$ and define the maps $\gammap$ and $\gammam$, as
above for $\ph$, to measure deviation from commutativity of the
diagrams
\[
\xymatrix@!C=4em@L=2ex{%
\sbar r\sbar s%
\ar[r]^{\sigmam(r,s)}%
\ar[d]_{\shat r\shat s}%
\ar@{}[dr]|{\gammam(r,s)}%
&\ol{rs}%
\ar[d]^{\wh{rs}}%
\\
\stil r\stil s%
\ar[r]_{\sigmam'(r,s)}%
&\wt{rs}%
}%
\]
and
\[
\xymatrix@!C=4em@L=2ex{%
\sbar r_0+\sbar r_1%
\ar[r]^{\sigmap(r_0,r_1)}%
\ar[d]_{\shat r_0+\shat r_1}%
\ar@{}[dr]|{\gammap(r_0,r_1)}%
&\ol{r_0+r_1}\ar[d]^{\wh{r_0+r_1}}%
\\
\stil r_0+\stil r_1%
\ar[r]_{\sigmap'(r_0,r_1)}%
&\wt{r_0+r_1}.%
}%
\]
That is, we define
\begin{align*}
\gammam(r,s)&=\wh{rs}\circ\sigmam(r,s)-\sigmam'(r,s)\circ\shat r\shat s%
\intertext{and}%
\gammap(r_0,r_1)&=\wh{r_0+r_1}\circ\sigmap(r_0,r_1)-\sigmap'(r_0,r_1)\circ(\shat r_0+\shat r_1).%
\end{align*}

The above four equalities then follow from considering the following
four diagrams, in view, as before, of \eqref{holep} and
\eqref{holem}, where ``$\comm$'' marks the strictly commuting
quadrangles:

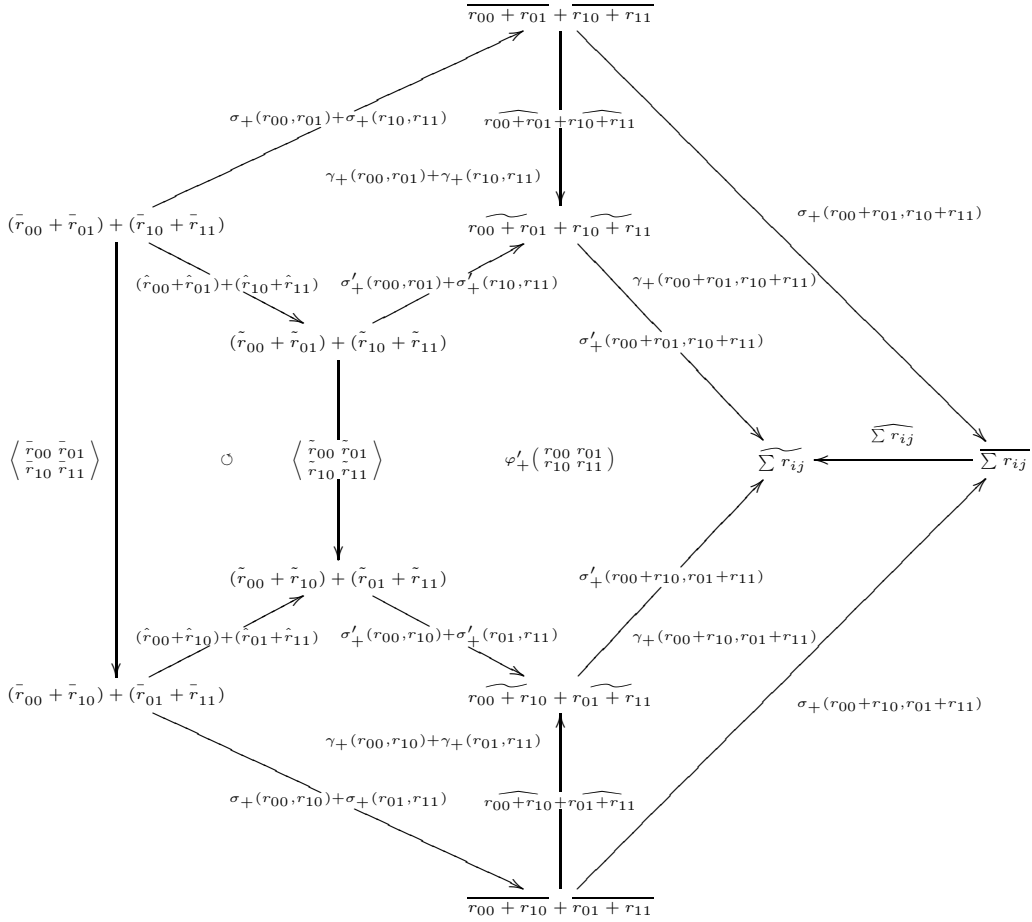
\begin{center}

\vfill

\scriptsize

\
\xymatrix{%
&&\sbar r(\sbar s\sbar t)%
\ar[dr]|{\shat r(\shat s\shat t)}%
\ar[rrrr]^{\sbar r\sigmam(s,t)}%
\ar@{}[dd]|\comm%
&&\ar@{}[d]|{r\gammam(s,t)}%
&&\sbar r\ol{st}%
\ar[ddrr]^{\sigmam(r,st)}%
\ar[dl]|{\shat r\wh{st}}%
\ar@{}[dd]|{\gammam(r,st)}%
\\
&&&\stil r(\stil s\stil t)%
\ar[rr]|{\stil r\sigmam'(s,t)}%
&&\stil r\wt{st}%
\ar[dr]|{\sigmam'(r,st)}%
\\
(\sbar r\sbar s)\sbar t%
\ar[uurr]^{\mass{\sbar r}{\sbar s}{\sbar t}}%
\ar[rr]|{(\shat r\shat s)\shat t}%
\ar[ddrrrr]_{\sigmam(r,s)\sbar t}%
&&(\stil r\stil s)\stil t%
\ar[ur]|{\mass{\stil r}{\stil s}{\stil t}}%
\ar[drr]|{\sigmam'(r,st)}%
\ar@{}[rrrr]^{\phm'(r,s,t)}%
\ar@{}[ddrr]|{\gammam(r,s)t}%
&&&&\wt{rst}%
\ar@{}[ddll]|{\gammam(rs,t)}%
&&\ol{rst}%
\ar[ll]|{\wh{rst}}%
\\
&&&&\wt{rs}\stil t%
\ar[urr]|{\sigmam'(rs,t)}
\\
&&&&\ol{rs}\sbar t%
\ar[u]|{\wh{rs}\shat t}%
\ar[uurrrr]_{\sigmam(rs,t)}%
}%
\hfill

\vfill

\qquad
\xymatrix@!C=2.3em@R=4em{%
&&\sbar r\sbar s_0+\sbar r\sbar s_1%
\ar[dr]|{\shat r\shat s_0+\shat r\shat s_1}%
\ar[rrrr]^{\sigmam(r,s_0)+\sigmam(r,s_1)}%
\ar@{}[dd]|{r\gammap(s_0,s_1)}%
&&\ar@{}[d]|{\gammam(r,s_0)+\gammam(r,s_1)}%
&&\ol{rs_0}+\ol{rs_1}%
\ar[ddrr]^{\sigmap(rs_0,rs_1)}%
\ar[dl]|{\wh{rs_0}+\wh{rs_1}}%
\ar@{}[dd]|{\gammap(rs_0,rs_1)}%
\\
&&&\stil r\stil s_0+\stil r\stil s_1%
\ar[rr]^{\sigmam'(r,s_0)+\sigmam'(r,s_1)}%
&&\wt{rs_0}+\wt{rs_1}%
\ar[dr]|{\sigmap'(rs_0,rs_1)}%
\\
\sbar r(\sbar s_0+\sbar s_1)%
\ar[uurr]^{\ldis{\sbar r}{\sbar s_0}{\sbar s_1}}%
\ar[rr]^{\shat r(\shat s_0+\shat s_1)}%
\ar[ddrrrr]_{\sbar r\sigmap(s_0,s_1)}%
&&\stil r(\stil s_0+\stil s_1)%
\ar[ur]|{\ldis{\stil r}{\stil s_0}{\stil s_1}}%
\ar[drr]|{r\sigmap'(s_0,s_1)}%
\ar@{}[rrrr]^{\ph_{\cdot+}'(r,s_0,s_1)}%
\ar@{}[ddrr]|\comm%
&&&&\wt{r(s_0+s_1)}%
\ar@{}[ddll]|{\gammam(r,s_0+s_1)}%
&&\ol{r(s_0+s_1)}%
\ar[ll]_{\wh{r(s_0+s_1)}}%
\\
&&&&\stil r\wt{s_0+s_1}%
\ar[urr]|{\sigmap'(r,s_0+s_1)}
\\
&&&&\sbar r\ol{s_0+s_1}%
\ar[u]|{\shat r\wh{s_0+s_1}}%
\ar[uurrrr]_{\sigmam(r,s_0+s_1)}%
}%
\hfill

\vfill

\qquad
\xymatrix@!C=2.3em@R=4em@L=2ex{%
&&\sbar r_0\sbar s+\sbar r_1\sbar s%
\ar[dr]|{\shat r_0\shat s+\shat r_1\shat s}%
\ar[rrrr]^{\sigmam(r_0,s)+\sigmam(r_1,s)}%
\ar@{}[dd]|{\gammap(r_0,r_1)s}%
&&\ar@{}[d]|{\gammam(r_0,s)+\gammam(r_1,s)}%
&&\ol{r_0s}+\ol{r_1s}%
\ar[ddrr]^{\sigmap(r_0s,r_1s)}%
\ar[dl]|{\wh{r_0s}+\wh{r_1s}}%
\ar@{}[dd]|{\gammap(r_0s,r_1s)}%
\\
&&&\stil r_0\stil s+\stil r_1\stil s%
\ar[rr]^{\sigmam'(r_0,s)+\sigmam'(r_1,s)}%
&&\wt{r_0s}+\wt{r_1s}%
\ar[dr]|{\sigmap'(r_0s,r_1s)}%
\\
(\sbar r_0+\sbar r_1)\sbar s%
\ar[uurr]^{\rdis{\sbar r_0}{\sbar r_1}{\sbar s}}%
\ar[rr]^{(\shat r_0+\shat r_1)\shat s}%
\ar[ddrrrr]_{\sigmap(r_0,r_1)\sbar s}%
&&(\stil r_0+\stil r_1)\stil s%
\ar[ur]|{\rdis{\stil r_0}{\stil r_1}{\stil s}}%
\ar[drr]|{\sigmap'(r_0,r_1)s}%
\ar@{}[rrrr]^{\ph_{+\cdot}'(r_0,r_1,s)}%
\ar@{}[ddrr]|\comm%
&&&&\wt{(r_0+r_1)s}%
\ar@{}[ddll]|{\gammam(r_0+r_1,s)}%
&&\ol{(r_0+r_1)s}%
\ar[ll]_{\wh{(r_0+r_1)s}}%
\\
&&&&\wt{r_0+r_1}\stil s%
\ar[urr]|{\sigmap'(r_0+r_1,s)}
\\
&&&&\ol{r_0+r_1}\sbar s%
\ar[u]|{\wh{r_0+r_1}\shat s}%
\ar[uurrrr]_{\sigmam(r_0+r_1,s)}%
}%
\hfill

\vfill

\tiny

\qquad
\xymatrix@!C=10em@R=5em@L=2ex{%
&&\ol{r_{00}+r_{01}}+\ol{r_{10}+r_{11}}%
\ar[dd]|{\wh{r_{00}+r_{01}}+\wh{r_{10}+r_{11}}}%
\ar[ddddrr]^{\sigmap(r_{00}+r_{01},r_{10}+r_{11})}%
\ar@{}[dddl]|{\gammap(r_{00},r_{01})+\gammap(r_{10},r_{11})\qquad}%
\ar@{}[ddddr]|(.6){\qquad\qquad_{\gammap(r_{00}+r_{01},r_{10}+r_{11})}}%
\\
\\
(\sbar r_{00}+\sbar r_{01})+(\sbar r_{10}+\sbar r_{11})%
\ar[uurr]|{\sigmap(r_{00},r_{01})+\sigmap(r_{10},r_{11})}%
\ar[dr]|{(\shat r_{00}+\shat r_{01})+(\shat r_{10}+\shat r_{11})}%
\ar[dddd]_{\coma{\sbar r_{00}}{\sbar r_{01}}{\sbar r_{10}}{\sbar r_{11}}}%
&&\wt{r_{00}+r_{01}}+\wt{r_{10}+r_{11}}%
\ar[ddr]|{\sigmap'(r_{00}+r_{01},r_{10}+r_{11})}%
\\
&(\stil r_{00}+\stil r_{01})+(\stil r_{10}+\stil r_{11})%
\ar[ur]|{\sigmap'(r_{00},r_{01})+\sigmap'(r_{10},r_{11})}%
\ar[dd]|{\coma{\stil r_{00}}{\stil r_{01}}{\stil r_{10}}{\stil r_{11}}}%
\\
\ar@{}[r]|\comm%
&\ar@{}[rr]|{\php'\left(\smat{r_{00}&r_{01}\\r_{10}&r_{11}}\right)}%
&&\wt{\sum r_{ij}}%
&\ol{\sum r_{ij}}%
\ar[l]_{\wh{\sum r_{ij}}}%
\\
&(\stil r_{00}+\stil r_{10})+(\stil r_{01}+\stil r_{11})%
\ar[dr]|{\sigmap'(r_{00},r_{10})+\sigmap'(r_{01},r_{11})}%
\\
(\sbar r_{00}+\sbar r_{10})+(\sbar r_{01}+\sbar r_{11})%
\ar[ddrr]|{\sigmap(r_{00},r_{10})+\sigmap(r_{01},r_{11})}%
\ar[ur]|{(\shat r_{00}+\shat r_{10})+(\shat r_{01}+\shat r_{11})}%
&&\wt{r_{00}+r_{10}}+\wt{r_{01}+r_{11}}%
\ar[uur]|{\sigmap'(r_{00}+r_{10},r_{01}+r_{11})}%
\\
\\
&&\ol{r_{00}+r_{10}}+\ol{r_{01}+r_{11}}%
\ar[uu]|{\wh{r_{00}+r_{10}}+\wh{r_{01}+r_{11}}}%
\ar[uuuurr]_{\sigmap(r_{00}+r_{10},r_{01}+r_{11})}%
\ar@{}[uuul]|{\gammap(r_{00},r_{10})+\gammap(r_{01},r_{11})\qquad}%
\ar@{}[uuuur]|(.6){\qquad\qquad_{\gammap(r_{00}+r_{10},r_{01}+r_{11})}}%
}%
\hfill

\end{center}

\end{proof}

We have thus obtained a map from the set of all categorical rings
$\sR$ with $\pi_0(\sR)=R$, $\pi_1(\sR)=B$ and matching bimodule
structure to the group $H^3(R;B)$. Let us next show that this map
factors through a quotient of the former set to yield a map
\[
\brk-:\mathrm{Crext}(R;B)\to H^3(R;B),
\]
where Crext$(R;B)$ denotes the set of equivalence classes of categorical
rings with $\pi_0$ equal to $R$ and $\pi_1$ equal to $B$, two such being
considered equivalent if there exists a 2-homomorphism between them
inducing identities on $R$ and $B$.

Indeed, in the same way as in \eqref{propcorr} we more generally
have:

\begin{Proposition}
For any categorical rings $\sR$ and $\sR'$ such that there exists a
2-homomorphism $\sR\to\sR'$ inducing identity maps on $\pi_0$ and $\pi_1$,
one has $\brk\sR=\brk{\sR'}$.
\end{Proposition}

\begin{proof}
Given a 2-homomorphism $\fb:\sR\to\sR'$, let us choose $\bar\ $ and
$\sigmam$, $\sigmap$ for $\sR$ as above, and then choose the
corresponding maps for $\sR'$ as follows:
\begin{align*}
\sbar r'&=f(\sbar r),\\
\sigmam'(r,s)&=\left(f\sbar rf\sbar s\xto{\fm(\sbar r,\sbar s)}f(\sbar r\sbar s)\xto{f\sigmam(r,s)}f\ol{rs}\right),\\
\sigmap'(r_0,r_1)&=\left(f\sbar r_0+f\sbar r_1\xto{\fp(\sbar r_0,\sbar
r_1)}f(\sbar r_0+\sbar r_1)\xto{f\sigmap(r_0,r_1)}f\ol{r_0+r_1}\right).
\end{align*}

In view of \eqref{holef}, \eqref{holep}, since $\fb$ induces
identity on $\pi_1$, i.~e. $f_\#$ is the identity map, one has the
diagrams

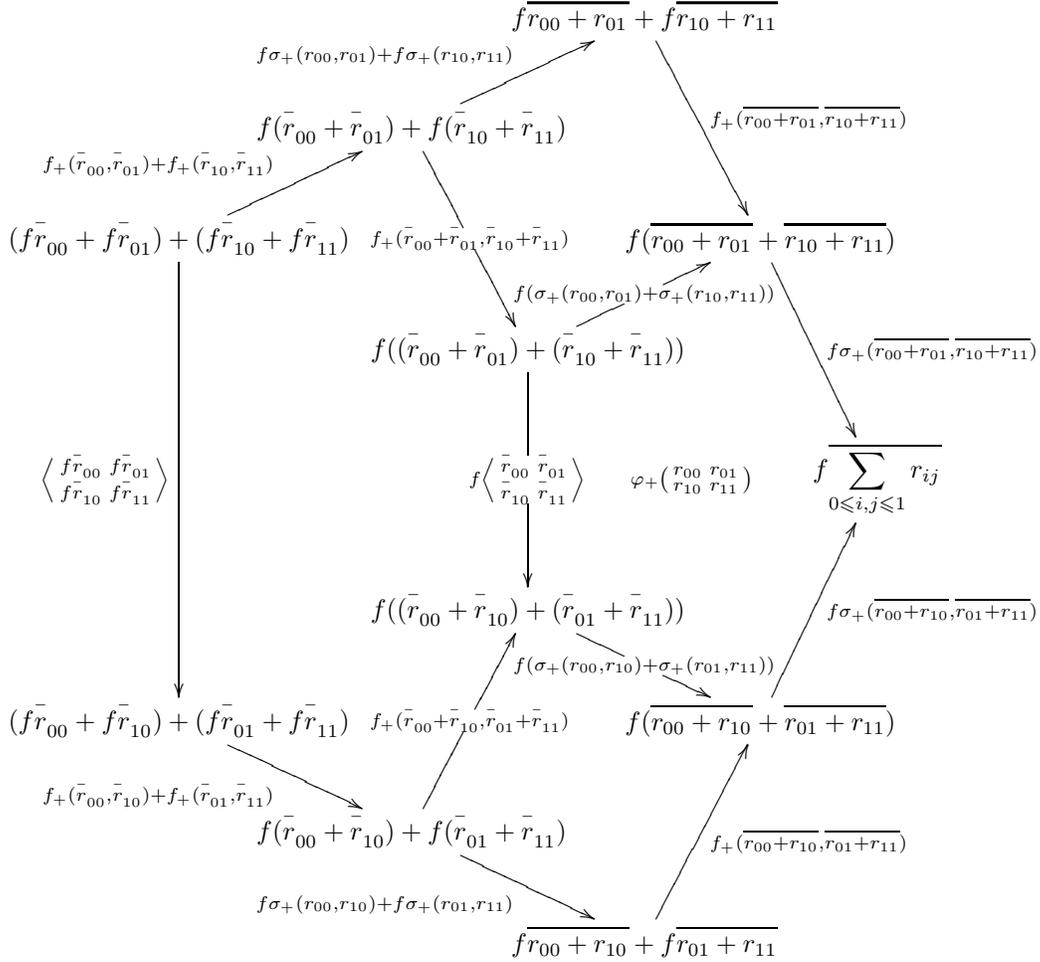
\[
\xymatrix@!C=2em{%
&&f\sbar r(f\sbar sf\sbar t)%
\ar[rr]^{f\sbar r\fm(\sbar s,\sbar t)}%
&&f\sbar rf(\sbar s\sbar t)%
\ar[dr]_{\fm(\sbar r,\sbar s\sbar t)}%
\ar[rr]^{f\sbar rf\sigmam(s,t)}%
&&f\sbar rf\ol{st}%
\ar[dr]^{\fm(\sbar r,\ol{st})}%
\\
&&&&&f(\sbar r(\sbar s\sbar t))%
\ar[rr]^{f(\sbar r\sigmam(s,t))}%
&&f(\sbar r\ol{st})%
\ar[dr]^{f\sigmam(r,st)}%
\\
(f\sbar rf\sbar s)f\sbar t%
\ar[uurr]^{\mass{f\sbar r}{f\sbar s}{f\sbar t}}%
\ar[drr]_{\fm(\sbar r,\sbar s)f\sbar t\qquad}%
&&&&f((\sbar r\sbar s)\sbar t)%
\ar[ur]^{f\mass{\sbar r}{\sbar s}{\sbar t}}%
\ar[drr]_{f(\sigmam(r,s)\sbar t)\quad}%
\ar@{}[rrrr]^{\phm(r,s,t)}%
&&&&f(\ol{rst})%
\\
&&f(\sbar r\sbar s)f\sbar t%
\ar[urr]_{\fm(\sbar r\sbar s,\sbar t)}%
\ar[drr]_{f\sigmam(r,s)f\sbar t\qquad}%
&&&&f(\ol{rs}\sbar t)%
\ar[urr]_{f\sigmam(rs,t)}%
\\
&&&&f\ol{rs}f\sbar t%
\ar[urr]_{\fm(\ol{rs},\sbar t)}%
}%
\]

\[
\xymatrix@!C=1.9em@L=2ex{%
&&f\sbar rf\sbar s_0+f\sbar rf\sbar s_1%
\ar[rr]^{\fm(\sbar r,\sbar s_0)+\fm(\sbar r,\sbar s_1)}%
&&f(\sbar r\sbar s_0)+f(\sbar r\sbar s_1)%
\ar[dr]|{\fp(\sbar r\sbar s_0,\sbar r\sbar s_1)}%
\ar[rr]^{f\sigmam(r,s_0)+f\sigmam(r,s_1)}%
&&f\ol{rs_0}+f\ol{rs_1}%
\ar[dr]^{\fp(\ol{rs_0},\ol{rs_1})}%
\\
&&&&&f(\sbar r\sbar s_0+\sbar r\sbar s_1)%
\ar[rr]^{f(\sigmam(r,s_0)+\sigmam(r,s_1))\qquad}%
&&f(\ol{rs_0}+\ol{rs_1})%
\ar[dr]^{f\sigmap(rs_0,rs_1)}%
\\
f\sbar r(f\sbar s_0+f\sbar s_1)%
\ar[uurr]^{\ldis{f\sbar r}{f\sbar s_0}{f\sbar s_1}}%
\ar[drr]_{f\sbar r\fp(\sbar s_0,\sbar s_1)}%
&&&&f(\sbar r(\sbar s_0+\sbar s_1))%
\ar[ur]^{f\ldis{\sbar r}{\sbar s_0}{\sbar s_1}}%
\ar[drr]|{f(\sbar r\sigmap(s_0,s_1))}%
\ar@{}[rrrr]^{\ph_{\cdot+}(r,s_0,s_1)}%
&&&&f(\ol{r(s_0+s_1)})%
\\
&&f\sbar rf(\sbar s_0+\sbar s_1)%
\ar[urr]|{\fm(\sbar r,\sbar s_0+\sbar s_1)}%
\ar[drr]_{f\sbar rf\sigmap(s_0,s_1)}%
&&&&f(\sbar r\ol{s_0+s_1})%
\ar[urr]_{f\sigmam(r,s_0+s_1)}%
\\
&&&&f\sbar rf\ol{s_0+s_1}%
\ar[urr]_{\fm(\sbar r,\ol{s_0+s_1})}%
}%
\]

\[
\xymatrix@!C=1.9em@L=2ex{%
&&f\sbar r_0f\sbar s+f\sbar r_1f\sbar s%
\ar[rr]^{\fm(\sbar r_0,\sbar s)+\fm(\sbar r_1,\sbar s)}%
&&f(\sbar r_0\sbar s)+f(\sbar r_1\sbar s)%
\ar[dr]|{\fp(\sbar r_0\sbar s,\sbar r_1\sbar s)}%
\ar[rr]^{f\sigmam(r_0,s)+f\sigmam(r_1,s)}%
&&f\ol{r_0s}+f\ol{r_1s}%
\ar[dr]^{\fp(\ol{r_0s},\ol{r_1s})}%
\\
&&&&&f(\sbar r_0\sbar s+\sbar r_1\sbar s)%
\ar[rr]^{f(\sigmam(r_0,s)+\sigmam(r_1,s))\qquad}%
&&f(\ol{r_0s}+\ol{r_1s})%
\ar[dr]^{f\sigmap(r_0s,r_1s)}%
\\
(f\sbar r_0+f\sbar r_1)f\sbar s%
\ar[uurr]^{\rdis{f\sbar r_0}{f\sbar r_1}{f\sbar s}}%
\ar[drr]_{\fp(\sbar r_0,\sbar r_1)f\sbar s}%
&&&&f((\sbar r_0+\sbar r_1)\sbar s)%
\ar[ur]^{f\rdis{\sbar r_0}{\sbar r_1}{\sbar s}}%
\ar[drr]|{f(\sigmap(r_0,r_1)\sbar s)}%
\ar@{}[rrrr]^{\ph_{+\cdot}(r_0,r_1,s)}%
&&&&f(\ol{(r_0+r_1)s})%
\\
&&f(\sbar r_0+\sbar r_1)f\sbar s%
\ar[urr]|{\fm(\sbar r_0+\sbar r_1,\sbar s)}%
\ar[drr]_{f\sigmap(r_0,r_1)f\sbar s}%
&&&&f(\ol{r_0+r_1}\sbar s)%
\ar[urr]_{f\sigmam(r_0+r_1,s)}%
\\
&&&&f\ol{r_0+r_1}f\sbar s%
\ar[urr]_{\fm(\ol{r_0+r_1},\sbar s)}%
}%
\]

\[
\xymatrix@!C=2em{%
&&&&f\ol{r_{00}+r_{01}}+f\ol{r_{10}+r_{11}}%
\ar[ddr]^{\fp(\ol{r_{00}+r_{01}},\ol{r_{10}+r_{11}})}%
\\
&&f(\sbar r_{00}+\sbar r_{01})+f(\sbar r_{10}+\sbar r_{11})%
\ar[urr]^{f\sigmap(r_{00},r_{01})+f\sigmap(r_{10},r_{11})\qquad\qquad}%
\ar[ddr]|{\fp(\sbar r_{00}+\sbar r_{01},\sbar r_{10}+\sbar r_{11})}%
\\
(f\sbar r_{00}+f\sbar r_{01})+(f\sbar r_{10}+f\sbar r_{11})%
\ar[urr]^{\fp(\sbar r_{00},\sbar r_{01})+\fp(\sbar r_{10},\sbar r_{11})\qquad\qquad}%
\ar[dddd]_{\coma{f\sbar r_{00}}{f\sbar r_{01}}{f\sbar r_{10}}{f\sbar r_{11}}}%
&&&&&f(\ol{r_{00}+r_{01}}+\ol{r_{10}+r_{11}})%
\ar[ddr]^{f\sigmap(\ol{r_{00}+r_{01}},\ol{r_{10}+r_{11}})}%
\ar@{}[dddd]_{\php\left(\smat{r_{00}&r_{01}\\r_{10}&r_{11}}\right)}%
\\
&&&f((\sbar r_{00}+\sbar r_{01})+(\sbar r_{10}+\sbar r_{11}))%
\ar[urr]|{f(\sigmap(r_{00},r_{01})+\sigmap(r_{10},r_{11}))}%
\ar[dd]|{f\coma{\sbar r_{00}}{\sbar r_{01}}{\sbar r_{10}}{\sbar r_{11}}}%
\\
&&&&&&f\ol{\displaystyle\sum_{0\le i,j\le1}r_{ij}}%
\\
&&&f((\sbar r_{00}+\sbar r_{10})+(\sbar r_{01}+\sbar r_{11}))%
\ar[drr]|{f(\sigmap(r_{00},r_{10})+\sigmap(r_{01},r_{11}))}%
\\
(f\sbar r_{00}+f\sbar r_{10})+(f\sbar r_{01}+f\sbar r_{11})%
\ar[drr]_{\fp(\sbar r_{00},\sbar r_{10})+\fp(\sbar r_{01},\sbar r_{11})\qquad\qquad}%
&&&&&f(\ol{r_{00}+r_{10}}+\ol{r_{01}+r_{11}})%
\ar[uur]_{f\sigmap(\ol{r_{00}+r_{10}},\ol{r_{01}+r_{11}})}%
\\
&&f(\sbar r_{00}+\sbar r_{10})+f(\sbar r_{01}+\sbar r_{11})%
\ar[drr]_{f\sigmap(r_{00},r_{10})+f\sigmap(r_{01},r_{11})\qquad\qquad}%
\ar[uur]|{\fp(\sbar r_{00}+\sbar r_{10},\sbar r_{01}+\sbar r_{11})}%
\\
&&&&f\ol{r_{00}+r_{10}}+f\ol{r_{01}+r_{11}}%
\ar[uur]_{\fp(\ol{r_{00}+r_{10}},\ol{r_{01}+r_{11}})}%
}%
\]
where unlabeled polygons strictly commute by coherence of $\fb$ and
naturality. These diagrams show that the cocycles $\ph$ and $\ph'$
representing characteristic classes of, respectively, $\sR$ and
$\sR'$, are cohomologous.
\end{proof}

We have thus defined a map
\[
\brk-:\mathrm{Crext}(R;B)\to H^3(R;B).
\]
Let us now construct a map in the opposite direction.

For a 3-cocycle $\ph=(\phm,\ph_{\cdot+},\ph_{+\cdot},\php)$ of $R$
with coefficients in $B$ let $\sR_\ph$ be the following categorical
ring. The set of objects of $\sR_\ph$ is $R$. The set of morphisms
is $B\x R$, where $(b,r)$ is a morphism from $r$ to $r$. Identities
are morphisms of the form $(0,r)$, and composition is given by
$(b,r)\circ(b',r)=(b+b',r)$. Categorical group structure is as
follows. Addition of objects and morphisms is given by
\[
(b_0,r_0)+(b_1,r_1)=(b_0+b_1,r_0+r_1);
\]
the neutral object is $(0,0)$, with the neutrality constraints given
by $\lambda(r)=\rho(r)=(0,r)$, the associativity constraint is given
by
\[
\aass rst=\left(\php\left(\smat{r&s\\0&t}\right),r+s+t\right),
\]
and the symmetry by
\[
\acom rs=\left(\php\left(\smat{0&r\\s&0}\right),r+s\right).
\]
Note that the latter two equalities are equivalent to the equality
\[
\coma{r_{00}}{r_{01}}{r_{10}}{r_{11}}=\left(\php\left(\smat{r_{00}&r_{01}\\r_{10}&r_{11}}\right),\sum_{0\le
i,j\le1}r_{ij}\right).
\]

Next, we define multiplication of objects and morphisms by
\[
(b,r)(b',r')=(br'+rb',rr'),
\]
the unit by $1$, the associativity constraint for the multiplication
by
\[
\mass rst=(\phm(r,s,t),rst)
\]
and the unitality constraints by
\begin{align*}
\lambda.(r)&=(-\phm(1,1,r),r),\\
\rho.(r)&=(\phm(r,1,1),r).
\end{align*}
Moreover, we define the distributivity constraints by the equalities
\[
\ldis r{s_0}{s_1}=(\ph_{\cdot+}(r,s_0,s_1),r(s_0+s_1))
\]
and
\[
\rdis {r_0}{r_1}s=(\ph_{+\cdot}(r_0,r_1,s),(r_0+r_1)s).
\]
It then turns out that commutativity of the coherence diagrams
necessary for $\sR_\ph$ to be a categorical ring correspond
precisely to the equations expressing the cocycle condition for
$\ph$.

Now suppose we are given two cohomologous 3-cocycles $\ph$, $\ph'$,
i.~e. there is a 2-cochain $\gamma=(\gammam,\gammap)$ satisfying the
required equalities. We then define a 2-homomorphism
$\fb=(f,\fp,\fm,f_1):\sR_\ph\to\sR_{\ph'}$ as follows. Since the
underlying categories of $\sR_\ph$ and $\sR_{\ph'}$ are identical,
we can define $f$ to be the identity functor. Moreover we define
\begin{align*}
\fp(r_0,r_1)&=(\gammap(r_0,r_1),r_0+r_1),\\
\fm(r,s)&=(\gammam(r,s),rs)%
\intertext{and}%
f_1&=(\gammam(1,1),1).
\end{align*}
Then again it is straightforward to verify that the coherence
conditions for $\fb$ to be a 2-homomorphism precisely amount to the
equalities expressing the fact that $\ph'$ differs from $\ph$ by the
coboundary of $\gamma$.

We have thus obtained a well-defined map
\[
\sR_-:H^3(R;B)\to\mathrm{Crext}(R;B)
\]
in the opposite direction.

Now it is obvious that constructing the characteristic class
$\brk{\sR_\ph}$ we can choose the maps $\sigmam$, $\sigmap$ in the
beginning of section \ref{ch} to be identities, which will produce
the cocycle $\ph$ back. Thus $\brk{\sR_\ph}$ is equal to the
cohomology class of $\ph$. So one composite of our maps (from $H^3$
to itself) is in fact identity. For the other composite to be also
the identity, it thus remains to construct, for any categorical ring
$\sR$, a 2-homomorphism between $\sR$ and $\sR_\ph$ for some cocycle
$\ph$, inducing identity on $\pi_0$ and $\pi_1$. For this, let us
return to the construction of the characteristic class of $\sR$; for
that construction, we have chosen an object $\bar r$ of $\sR$ in
each isomorphism class $r\in\pi_0(\sR)=R$ and morphisms $\sigmam$,
$\sigmap$, which then produced the cocycle $\ph$ representing
$\brk{\sR}$. Obviously these choices can be made in such a way that
$\bar0=0$, $\bar1=1$, $\sigmap(0,r)=\lambda(\bar r)$,
$\sigmap(r,0)=\rho(\bar r)$, $\sigmam(1,r)=\lambda.(\bar r)$, and
$\sigmam(r,1)=\rho.(\bar r)$. Let us then use these data to define a
functor $f:\sR_\ph\to\sR$. On objects this functor is given by
$f(r)=\bar r$ and on morphisms by
\[
f(b,r)=\bar r\xto{\lambda(\bar r)\1}0+\bar r\xto{b+\bar r}0+\bar
r\xto{\lambda(\bar r)}\bar r.
\]
This $f$ then extends to a 2-homomorphism
$\fb=(f,\fp,\fm,f_1):\sR_\ph\to\sR$ with $\fm=\sigmam$,
$\fp=\sigmap$ and $f_1=$ identity of 1.

Summarizing all of the above, we have thus proved

\begin{Theorem}\label{main}
For any ring $R$ and any $R$-bimodule $B$ there is a bijection
\[
H^3(R;B)\approx\mathrm{Crext}(R;B)
\]
between the third Mac Lane cohomology of $R$ with coefficients in
$B$ and equivalence classes of categorical rings $\sR$ with
$\pi_0(\sR)=R$, $\pi_1(\sR)=B$ and the resulting bimodule structure
coinciding with the original one.
\end{Theorem}\qed

\end{document}